\documentclass[12pt,fullpage]{article}

\usepackage{amsfonts,xr,graphicx,amsmath,amssymb,colortbl,epsfig,subfigure}
\def\Z{\mathbb{Z}}
\def\R{\mathbb{R}}
\def\N{\mathbb{N}}

\def\tilde{\widetilde}
\def\epsilon{\varepsilon}
   \def\ga{\gamma}

\def\ee{\epsilon}
\def\beq{\begin{equation}}
\def\eeq{\end{equation}}

\def\div{\rm{div}}

\def\osc{\rm{osc}}
\def\length{\rm{length}}

\newcommand{\SE}{\setcounter{equation}{0} \section}
\newcommand{\baa}{\begin{array}}
\newcommand{\eaa}{\end{array}}
\newcommand{\ba}{\begin{eqnarray}}
\newcommand{\ea}{\end{eqnarray}}
\newtheorem{theo}{\bf Theorem}[section]
\newtheorem{lem}[theo]{\bf Lemma}
\newtheorem{pro}[theo]{\bf Proposition}
\newtheorem{cor}[theo]{\bf Corollary}

\newtheorem{rem}[theo]{\bf Remark}

\begin{document}

\title{\bf{A Liouville theorem for the Euler equations in the plane}}
\author{Fran{\c{c}}ois Hamel and Nikolai Nadirashvili$\,$\thanks{This work has been carried out in the framework of Archim\`ede Labex (ANR-11-LABX-0033) and of the A*MIDEX project (ANR-11-IDEX-0001-02), funded by the ``Investissements d'Avenir" French Government program managed by the French National Research Agency (ANR). The research leading to these results has also received funding from the European Research Council under the European Union's Seventh Framework Programme (FP/2007-2013) ERC Grant Agreement n.~321186~- ReaDi~- Reaction-Diffusion Equations, Propagation and Modelling and from the ANR NONLOCAL project (ANR-14-CE25-0013).}\\
\\
\footnotesize{Aix Marseille Univ, CNRS, Centrale Marseille, I2M, Marseille, France}}
\date{}
\maketitle

\begin{abstract}
This paper is concerned with qualitative properties of bounded steady flows of an ideal incompressible fluid with no stagnation point in the two-dimensional plane $\R^2$. We show that any such flow is a shear flow, that is, it is parallel to some constant vector. The proof of this Liouville-type result is firstly based on the study of the geometric properties of the level curves of the stream function and secondly on the derivation of some estimates on the at-most-logarithmic growth of the argument of the flow in large balls. These estimates lead to the conclusion that the streamlines of the flow are all parallel lines.
\end{abstract}



\SE{Introduction and main results}\label{intro}

In this paper, we consider steady flows $v=(v_1,v_2)$ of an ideal fluid in the two-dimensional plane~$\R^2$, which solve the system of the incompressible Euler equations: 
\begin{equation}\label{1}
\left\{\begin{array}{ll}
v\cdot\nabla\,v +\nabla\,p =0 & \mbox{in $\R^2$},\vspace{3pt}\\
{\div}\ v=0 & \mbox{in $\R^2$}.
\end{array} \right.
\end{equation}
Throughout the paper, the solutions are always understood in the classical sense, that is, $v$ and~$p$ are (at least) of class $C^1(\R^2)$ and satisfy~\eqref{1}. Any flow $v$ is called a shear flow if there is a unit vector $e=(e_1,e_2)\in\mathbb{S}^1$ such that $v$ is parallel to $e$ ($\mathbb{S}^1$ denotes the unit circle in $\R^2$). Due to the incompressibility condition ${\div}\,v=0$, any shear flow parallel to $e$ only depends on the orthogonal variable $x\cdot e^{\perp}$, where $e^{\perp}=(-e_2,e_1)$. In other words, a shear flow is a flow for which there are $e\in\mathbb{S}^1$ and a function $V:\R\to\R$ such that
\beq\label{shear}
v(x)=V(x\cdot e^{\perp})\,e
\eeq
for all $x\in\R^2$. It is easy to see that $v$ is a shear flow if and only if the pressure $p$ is constant. In the sequel, we denote $x\mapsto|x|$ the Euclidean norm in $\R^2$.\par
The main result of this paper is the following rigidity result for the stationary Euler equations.

\begin{theo}\label{th1}
Let $v$ be a $C^2(\R^2)$ flow solving~\eqref{1}. Assume that $v\in L^{\infty}(\R^2)$ and that
\beq\label{hyp1}
\inf_{\R^2}|v|>0.
\eeq
Then $v$ is a shear flow. Namely, $v$ is of the type~\eqref{shear} and the function $V$ in~\eqref{shear} has a constant strict sign.
\end{theo}

Theorem~\ref{th1} can also be viewed as a Liouville-type rigidity result since the conclusion says that the argument of the flow is actually constant, and that the pressure $p$ is constant as well.\par
Some comments on the assumptions made in Theorem~\ref{th1} are in order. First of all, the assumption~\eqref{hyp1} means that the flow $v$ has no stagnation point in $\R^2$ or at infinity. In other words, Theorem~\ref{th1} means that any $C^2(\R^2)\cap L^\infty(\R^2)$ flow which is not a shear flow must have a stagnation point in $\R^2$ or at infinity.\par
Without the condition~\eqref{hyp1}, the conclusion of Theorem~\ref{th1} does not hold in general. For instance, for any $(\alpha,\beta)\in\R^*\times\R^*$, the smooth cellular flow $v$ defined in $\R^2$ by
$$v(x_1,x_2)=\nabla^{\perp}\big(\sin(\alpha x_1)\sin(\beta x_2)\big)=\big(\!-\beta\sin(\alpha x_1)\cos(\beta x_2),\alpha\cos(\alpha x_1)\sin(\beta x_2)\big),$$
which solves~\eqref{1} with $p(x)=(\beta^2/4)\cos(2\alpha x_1)+(\alpha^2/4)\cos(2\beta x_2)$, is bounded, but it has (countably many) stagnation points in $\R^2$, and it is not a shear flow. However, we point out that, the sufficient condition~\eqref{hyp1} is obviously not equivalent to being a shear flow. Indeed, any continuous shear flow $v(x)=V(x\cdot e^{\perp})\,e$ for which $V$ changes sign (or more generally if $\inf_\R|V|=0$) does not satisfy the condition~\eqref{hyp1}.\par
Moreover, without the boundedness of $v$, the conclusion of Theorem~\ref{th1} does not hold either in general. For instance, the smooth flow $v$ defined in $\R^2$ by
$$v(x)=\nabla^{\perp}\big(x_2\cosh(x_1))=(-\cosh(x_1),x_2\sinh(x_1)),$$
which solves~\eqref{1} with $p(x)=-\cosh(2x_1)/4+x_2^2/2$, satisfies $\inf_{\R^2}|v|>0$ but it is not bounded in~$\R^2$, and it is not a shear flow.\par
The assumption on the $C^2(\R^2)$ smoothness of $v$ is a technical assumption which is used in the proof. It is connected with the $C^1$ smoothness of the nonlinear source term $f$ in the equation satisfied by the stream function $u$ of the flow $v$, see~\eqref{defu} and~\eqref{eqelliptic} below. We refer to Section~\ref{sec2} for further details. However, no uniform smoothness is assumed, namely $v$ is not assumed to be uniformly continuous and its first and second order derivatives are not assumed to be bounded nor uniformly continuous.\par
In our previous paper~\cite{hn}, we considered the case of a two-dimensional strip with bounded section and the case of the half-plane, assuming in both cases that the flow was tangential on the boundary. In those both situations, the boundary of the domain was a streamline  and the conclusion was that the flow is a shear flow, parallel to the boundary of the domain. In the present paper, there is no boundary and no obvious simple streamline. We shall circumvent this difficulty by proving additional estimates on the flow and its stream function at infinity, and in particular on the at-most-logarithmic growth of the argument of the flow in large balls. We also mention other rigidity results for the stationary solutions of~\eqref{1} in other two-dimensional domains, such as the analyticity of the streamlines under a condition of the type $v_1>0$ in the unit disc~\cite{kn}, and the local correspondence between the vorticities of the stationary solutions of~\eqref{1} and the co-adjoint orbits of the vorticities for the non-stationary version of~\eqref{1} in annular domains~\cite{cs}.\par
Lastly, to complete Section~\ref{intro}, we point out two immediate corollaries of Theorem~\ref{th1}. The first one is concerned with periodic flows. In the sequel, we say that a flow $v$ is periodic if there is a basis $(\mathrm{e}_1,\mathrm{e}_2)$ of $\R^2$ such that $v(x)=v(x+k_1\mathrm{e}_1+k_2\mathrm{e}_2)$ in $\R^2$ for all $x\in\R^2$ and $(k_1,k_2)\in\Z^2$.

\begin{cor}
Let $v$ be a $C^2(\R^2)$ periodic flow solving~\eqref{1}. If $|v(x)|\neq0$ for all $x\in\R^2$, then $v$ is a shear flow.
\end{cor}

The second corollary states that the class of bounded shear flows satisfying~\eqref{hyp1} is stable under small $L^{\infty}(\R^2)$ perturbations.

\begin{cor}
Let $v$ be a bounded shear flow solving~\eqref{1} and satisfying~\eqref{hyp1}. There is $\varepsilon>0$ such that, if $v'$ is a $C^2(\R^2)$ flow solving~\eqref{1} and satisfying $\|v'-v\|_{L^{\infty}(\R^2)}\le\varepsilon$, then $v'$ is a shear flow.
\end{cor}

\begin{rem}\label{remmonotone}{\rm If, in addition to the condition $\inf_{\R^2}|v|>0$, one assumes that $v\cdot e>0$ in $\R^2$ for some direction $e\in\mathbb{S}^1$ (by continuity, up to changing $v$ into $-v$, it is therefore sufficient to assume that $v(x)\cdot e\neq 0$ for all $x\in\R^2$), then the end of the proof of Theorem~\ref{th1} would be much simpler: indeed, in that case, the stream function $u$ defined in~\eqref{defu} below would be monotone in the direction $e^{\perp}$. Since $u$ satisfies a semilinear elliptic equation of the type $\Delta u+f(u)=0$ in $\R^2$ (see~\eqref{eqelliptic} below), it would then follow that $u$ is one-dimensional, as in the proof of a related conjecture of De Giorgi~\cite{dg} in dimension~$2$ (see~\cite{bcn,gg} and see also~\cite{aac,ac,dkw,f2,fv,s} for further references in that direction). Finally, since $u$ is one-dimensional, the vector field $v$ is a shear flow. We refer to Section~\ref{sec23} below for further details.}
\end{rem}

\begin{rem} {\rm A Liouville theorem is known for the Navier-Stokes equations on the plane~\cite{knss}. For a viscous flow the Liouville property has a different form: any uniformly bounded solution of the Navier-Stokes equations on the plane is a constant.}
\end{rem}


\SE{Proof of Theorem~\ref{th1}}\label{sec2}

Sections~\ref{sec20} and~\ref{sec21} are devoted to some important notations and to the proof of some preliminary lemmas. The proof of Theorem~\ref{th1} is completed in Section~\ref{sec22}, assuming the technical Proposition~\ref{growth} below on the at-most-logarithmic growth of the argument of the flow in large balls. In Section~\ref{sec23}, we consider the special case where $v\cdot e>0$ in $\R^2$ (see Remark~\ref{remmonotone} above), in which case the end of the proof of Theorem~\ref{th1} is much easier and does not require Proposition~\ref{growth}.


\subsection{The main scheme of the proof and some important notations}\label{sec20}

Let us first explain the main lines of the proof of Theorem~\ref{th1}. It is based on the study of the geometric properties of the streamlines of the flow $v$ and of the orthogonal trajectories of the gradient flow defined by the potential $u$ of the flow $v$ (see definition~\eqref{defu} below). The first main point is to show that all streamlines of $v$ are unbounded and foliate the plane $\R^2$ in a monotone way. Since the vorticity
$$\frac{\partial v_2}{\partial x_1}-\frac{\partial v_1}{\partial x_2}$$
is constant along the streamlines of the flow $v$, the potential function $u$ will be proved to satisfy a semilinear elliptic equation of the type $\Delta u+f(u)=0$ in $\R^2$. Another key-point consists in proving that the argument of the flow $v$ (and of $\nabla u$) grows at most as $\ln R$ in balls of large radius $R$. Finally, we use a compactness argument and a result of Moser~\cite{m} to conclude that the argument of $v$, which solves a uniformly elliptic linear equation in divergence form, is actually constant.\par
Throughout Section~\ref{sec2}, $v$ is a given $C^2(\R^2)\cap L^\infty(\R^2)$ vector field solving~\eqref{1}, and such that $\inf_{\R^2}|v|>0$. Therefore, there is $0<\eta\le1$ such that
\beq\label{defee0}
0<\eta\le|v(x)|\le\eta^{-1}\ \hbox{ for all }x\in\R^2.
\eeq
Our goal is to show that $v$ is a shear flow.\par
To do so, let us first introduce some important definitions. Let $u$ be a potential function (or stream function) of the flow $v$. More precisely, $u:\R^2\to\R$ is a $C^3(\R^2)$ function such that
\beq\label{defu}
\nabla^{\perp}u=v,\ \hbox{ that is},\ \ \frac{\partial u}{\partial x_1}=v_2\ \hbox{ and }\ \frac{\partial u}{\partial x_2}=-v_1
\eeq
in $\R^2$. Since $v$ is divergence free and $\R^2$ is simply connected, it follows that the potential function~$u$ is well and uniquely defined in $\R^2$ up to a constant. In the sequel, we call $u$ the unique stream function such that $u(0)=0$.\par
The trajectories of the flow $v$, that is, the curves tangent to $v$ at each point, are called the streamlines of the flow. Since $|v|>0$  in $\R^2$, a given streamline $\Gamma$ of $v$ cannot have an endpoint in~$\R^2$ and it always admits a $C^1$ parametrization $\gamma:\R\to\R^2$ ($\gamma(\R)=\Gamma$) such that $|\dot{\gamma}(t)|>0$ for all $t\in\R$. Actually, since $v$ is bounded too, it follows that, for any given $x\in\R^2$, the solution $\gamma_x$ of
\beq\label{defgammax}\left\{\baa{rcl}
\dot\gamma_x(t) & = & v(\gamma_x(t)),\vspace{3pt}\\
\gamma_x(0) & = & x\eaa\right.
\eeq
is defined in the whole interval $\R$ and is a parametrization of the streamline $\Gamma_x$ of $v$ containing~$x$. By definition, the stream function $u$ is constant along the streamlines of the flow $v$ and, for any given $x\in\R^2$, the level curve of $u$ containing $x$, namely the connected component of the level set~$\{y\in\R^2;\ u(y)=u(x)\}$ containing $x$, is nothing but the streamline $\Gamma_x$.\par
We will also consider in the proof of Theorem~\ref{th1} the trajectories of the gradient flow $\dot{\sigma}=\nabla u(\sigma)$. Namely, for any~$x\in\R^2$, let $\sigma_x$ be the solution of
\beq\label{xix}\left\{\baa{rcl}
\dot\sigma_x(t) & = & \nabla u(\sigma_x(t)),\vspace{3pt}\\
\sigma_x(0) & = & x.\eaa\right.
\eeq
As for $\gamma_x$ in~\eqref{defgammax}, the parametrization $\sigma_x$ of the trajectory $\Sigma_x$ of the gradient flow containing $x$ is defined in the whole $\R$. Furthermore, $\Sigma_x$ is orthogonal to $\Gamma_x$ at $x$.\par
In the sequel, we denote
$$B(x,r)=\big\{y\in\R^2;\ |x-y|<r\big\}$$
the open Euclidean ball of centre $x\in\R$ and radius $r>0$. We also use at some places the notation $0=(0,0)$ and then
$$B(0,r)=B((0,0),r)$$
for $r>0$.


\subsection{Some preliminary lemmas}\label{sec21}

Let us now establish a few fundamental elementary properties of the streamlines of the flow and of the trajectories of the gradient flow. The first such property is the unboundedness of the trajectories of the gradient flow.

\begin{lem}\label{unboundedbis}
Let $\Sigma$ be any trajectory of the gradient flow and let $\sigma:\R\to\R^2$ be any $C^1$ parametrization of~$\Sigma$ such that $|\dot{\sigma}(t)|>0$ for all $t\in\R$. Then $|\sigma(t)|\to+\infty$ as $|t|\to+\infty$ and the map $t\mapsto u(\sigma(t))$ is a homeomorphism from $\R$ to $\R$.
\end{lem}

\noindent{\bf{Proof.}} Let $x$ be any point on $\Sigma$, that is, $\Sigma=\Sigma_x$, and let $\sigma_x:\R\to\R^2$ be the parametrization of $\Sigma$ defined by~\eqref{xix}. The function $g:t\mapsto g(t):=u(\sigma_x(t))$ is (at least) of class $C^1(\R)$ and
\beq\label{increasingg}
g'(t)=|\nabla u(\sigma_x(t))|^2=|v(\sigma_x(t))|^2\ge\eta^2>0
\eeq
for all $t\in\R$ by~\eqref{defee0}. In particular, $g(t)\to\pm\infty$ as $t\to\pm\infty$. Since $u$ is locally bounded, it also follows that $|\sigma_x(t)|\to+\infty$ as $|t|\to+\infty$.\par
Now, consider any $C^1$ parametrization $\sigma:\R\to\R^2$ of~$\Sigma$ such that $|\dot{\sigma}(t)|>0$ for all $t\in\R$. Since $\dot{\sigma}(t)\neq(0,0)$ is parallel to $\nabla u(\sigma(t))\neq(0,0)$ at each $t\in\R$, the continuous function $t\mapsto\dot{\sigma}(t)\cdot\nabla u(\sigma(t))$ has a constant sign in $\R$ and the $C^1$ function $t\mapsto u(\sigma(t))$ is then either increasing or decreasing. Since $\sigma$ is a parametrization of $\Sigma$ and since $u(\Sigma)=\R$ from the previous paragraph, one infers that either $u(\sigma(t))\to\pm\infty$ as $t\to\pm\infty$, or $u(\sigma(t))\to\mp\infty$ as $t\to\pm\infty$. In any case, $t\mapsto u(\sigma(t))$ is a homeomorphism from $\R$ to $\R$ and, as for $\sigma_x(t)$, one gets that $|\sigma(t)|\to+\infty$ as $|t|\to+\infty$.\hfill$\Box$\break

An immediate consequence of the proof of Lemma~\ref{unboundedbis} is the following estimate, which we state separately since it will be used several times in the sequel.

\begin{lem}\label{lemsigma}
Let $\Sigma$ be any trajectory of the gradient flow, with parametrization $\sigma:\R\to\R^2$ solving $\dot\sigma(t)=\nabla u(\sigma(t))$ for all $t\in\R$. Then, $t\mapsto u(\sigma(t))$ is increasing and, for any real numbers $\alpha\neq\beta$, 
\beq\label{ualphabeta}
|u(\sigma(\alpha))-u(\sigma(\beta))|\ge\eta\times{\length}(\sigma(I))\ge\eta\,|\sigma(\alpha)-\sigma(\beta)|,
\eeq
where $I=\big[\min(\alpha,\beta),\max(\alpha,\beta)\big]$.
\end{lem}

\noindent{\bf{Proof.}} The fact that $g:t\mapsto u(\sigma(t))$ is increasing follows from~\eqref{increasingg}. Furthermore, for all $t\in\R$, $g'(t)=|\nabla u(\sigma(t))|^2\ge\eta\,|\dot\sigma(t)|$. This inequality immediately yields the conclusion.\hfill$\Box$\break

Like the trajectories of the gradient flow, it turns out that the streamlines of the flow are also unbounded. More precisely, the following result holds.

\begin{lem}\label{unbounded}
Let $\Gamma$ be any streamline of the flow $v$ and $\gamma:\R\to\R^2$ be any $C^1$ parametrization of~$\Gamma$ such that $|\dot{\gamma}(t)|>0$ for all $t\in\R$. Then $|\ga(t)|\to+\infty $ as $|t|\to+\infty$. 
\end{lem}

\noindent{\bf{Proof.}} First of all, let us notice that $\Gamma$ is not closed in the sense that any $C^1$ parametrization $\gamma:\R\to\R^2$ of $\Gamma$ such that $|\dot{\gamma}(t)|>0$ for all $t\in\R$ is actually one-to-one. Indeed, otherwise, there would exist two real numbers $a<b$ such that $\gamma(a)=\gamma(b)$ and $\Gamma$ would then be equal to $\gamma([a,b])$. Then the open set $\Omega$ surrounded by $\Gamma=\gamma([a,b])$ would be nonempty (since $|\dot{\gamma}|>0$ in $\R$) while, by definition of $u$, the function $u$ is constant on the streamline $\Gamma=\partial\Omega$. Thus, $u$ would have either an interior minimum or an interior maximum in $\Omega$, which is ruled out since $\nabla u=-v^{\perp}$ does not vanish.\par
Let now $x$ be any point on $\Gamma$ (in other words, $\Gamma=\Gamma_x$) and let $\gamma_x:\R\to\R^2$ be the parametrization of $\Gamma$ defined by~\eqref{defgammax}. We claim that $|\ga_x(t)|\to+\infty $ as $|t|\to+\infty$. Assume not. Then there are~$y\in\R^2$ and a sequence $(\tau_n)_{n\in\N}$ in $\R$ such that
$$|\tau_n|\to+\infty\ \hbox{ and }\ \gamma_x(\tau_n)\to y\ \hbox{ as }n\to+\infty.$$
Since $u(\gamma_x(\tau_n))=u(x)$ by definition of $\gamma_x$ and $u$, the continuity of $u$ implies that $u(x)=u(y)$. Call now
$$\Lambda=\big\{z\in\R^2;\ u(z)=u(x)=u(y)\big\}$$
the level set of $u$ with level $u(x)=u(y)$. Since $\nabla u(y)=-v^{\perp}(y)\neq(0,0)$, the implicit function theorem yields the existence of $r>0$ small enough such that $B(y,r)\cap\Lambda$ is a graph in the variable parallel to $v(y)$ and such that it can be written as
$$B(y,r)\cap\Lambda=B(y,r)\cap\Gamma_y=\big\{\gamma_y(t);\ \alpha<t<\beta\big\}$$
for some real numbers $\alpha<0<\beta$. Since $u(\gamma_x(\tau_n))=u(x)=u(y)$ and $\gamma_x(\tau_n)\to y$ as $n\to+\infty$, it follows that $\gamma_x(\tau_n)\in B(y,r)\cap\Lambda$ for $n$ large enough, hence $\gamma_x(\tau_n)\in\Gamma_y$ for such $n$. But $\gamma_x(\tau_n)\in\Gamma_x$ by definition, hence $\Gamma_x=\Gamma_y$ and
$$y=\gamma_x(t)$$
for some $t\in\R$. Furthermore, still for $n$ large enough, $\gamma_x(\tau_n)\in B(y,r)\cap\Lambda$, hence $\gamma_x(\tau_n)=\gamma_y(t_n)$ for some $t_n\in(\alpha,\beta)$. As a consequence, $\gamma_x(\tau_n)=\gamma_x(t+t_n)$ for $n$ large enough, and $\tau_n=t+t_n$ for such $n$ since $\gamma_x$ is one-to-one from the previous paragraph. This leads to a contradiction as $n\to+\infty$, since the sequence $(t_n)_{n\in\N}$ is bounded whereas $|\tau_n|\to+\infty$ as $n\to+\infty$.\par
Thus, $|\gamma_x(t)|\to+\infty$ as $|t|\to+\infty$. The same conclusion immediately follows for any $C^1$ parametrization $\gamma:\R\to\R^2$ of $\Gamma$ such that $|\dot{\gamma}(t)|>0$ for all $t\in\R$ and the proof of Lemma~\ref{unbounded} is thereby complete.\hfill$\Box$\break

The previous lemma says that the parametrizations $\gamma_x(t)$ solving $\dot\gamma_x(t)=v(\gamma_x(t))$ with $\gamma_x(0)=x$ converge to infinity in norm for any given $x\in\R^2$. It actually turns out that this convergence holds uniformly in $x$ when $x$ belongs to a fixed bounded set. Namely, the following result holds.

\begin{lem}\label{unboundedter}
Let $K\subset\R^2$ be a bounded set and let $\gamma_x:\R\to\R^2$ be defined as in~\eqref{defgammax}. Then $|\gamma_x(t)|\to+\infty$ as $|t|\to+\infty$ uniformly in $x\in K$.
\end{lem}

\noindent{\bf{Proof.}} Assume by way of contradiction that the conclusion does not hold for some bounded set $K\subset\R^2$. Then there are $x,y\in\R^2$, a sequence $(x_n)_{n\in\N}$ in $\R^2$ and a sequence $(t_n)_{n\in\N}$ in $\R$ such that
\beq\label{xnyn}
x_n\to x,\ \ y_n:=\gamma_{x_n}(t_n)\to y\ \hbox{ and }\ |t_n|\to+\infty\ \hbox{ as }n\to+\infty.
\eeq
Up to extraction of a subsequence, one can assume that $t_n\to\pm\infty$ as $n\to+\infty$. We only consider the case
\beq\label{tn}
t_n\to+\infty\ \hbox{ as }n\to+\infty
\eeq
(the case $t_n\to-\infty$ can be handled similarly).\par
Fix now some positive real numbers $r$ and then $R$ such that
\beq\label{rR}
\max\big(|x|,|y|\big)<r<R\ \hbox{ and }\ 2\,M_r<\eta\,(R-r),
\eeq
where $\eta>0$ is as in~\eqref{defee0} and
$$M_r:=\max_{\overline{B(0,r)}}|u|.$$
Then, since $|\gamma_y(t)|\to+\infty$ as $|t|\to+\infty$ by Lemma~\ref{unbounded}, there is $T>0$ such that $|\gamma_y(-T)|>R$.\par
From the continuous dependence of the solutions $\gamma_X(t)$ of~\eqref{defgammax} with respect to the initial value $X$ for each given $t\in\R$, one knows that $\gamma_{y_n}(-T)\to\gamma_y(-T)$ as $n\to+\infty$. Since $y_n=\gamma_{x_n}(t_n)$, this means that $\gamma_{x_n}(t_n-T)=\gamma_{y_n}(-T)\to\gamma_y(-T)$ as $n\to+\infty$, hence $|\gamma_{x_n}(t_n-T)|\to|\gamma_y(-T)|>R$ as $n\to+\infty$. Together with~\eqref{xnyn},~\eqref{tn} and~\eqref{rR}, one infers the existence of an integer $N\in\N$ such that
\beq\label{defn0}
|\gamma_{x_N}(0)|=|x_N|<r,\ \ |\gamma_{x_N}(t_N)|<r,\ \ 0<t_N-T<t_N\ \hbox{ and }\ |\gamma_{x_N}(t_N-T)|>R\ (>r).
\eeq
By continuity of the parametrization $t\mapsto\gamma_{x_N}(t)$ of $\Gamma_{x_N}$, there are some real number $\alpha$ and $\beta$ such that
$$0<\alpha<t_N-T<\beta<t_N,\ \ |\gamma_{x_N}(\alpha)|=|\gamma_{x_N}(\beta)|=r\ \hbox{ and }\ |\gamma_{x_N}(t)|>r\hbox{ for all }t\in(\alpha,\beta),$$
see Figure~1.\par
\begin{figure}\centering
\includegraphics[scale=0.8]{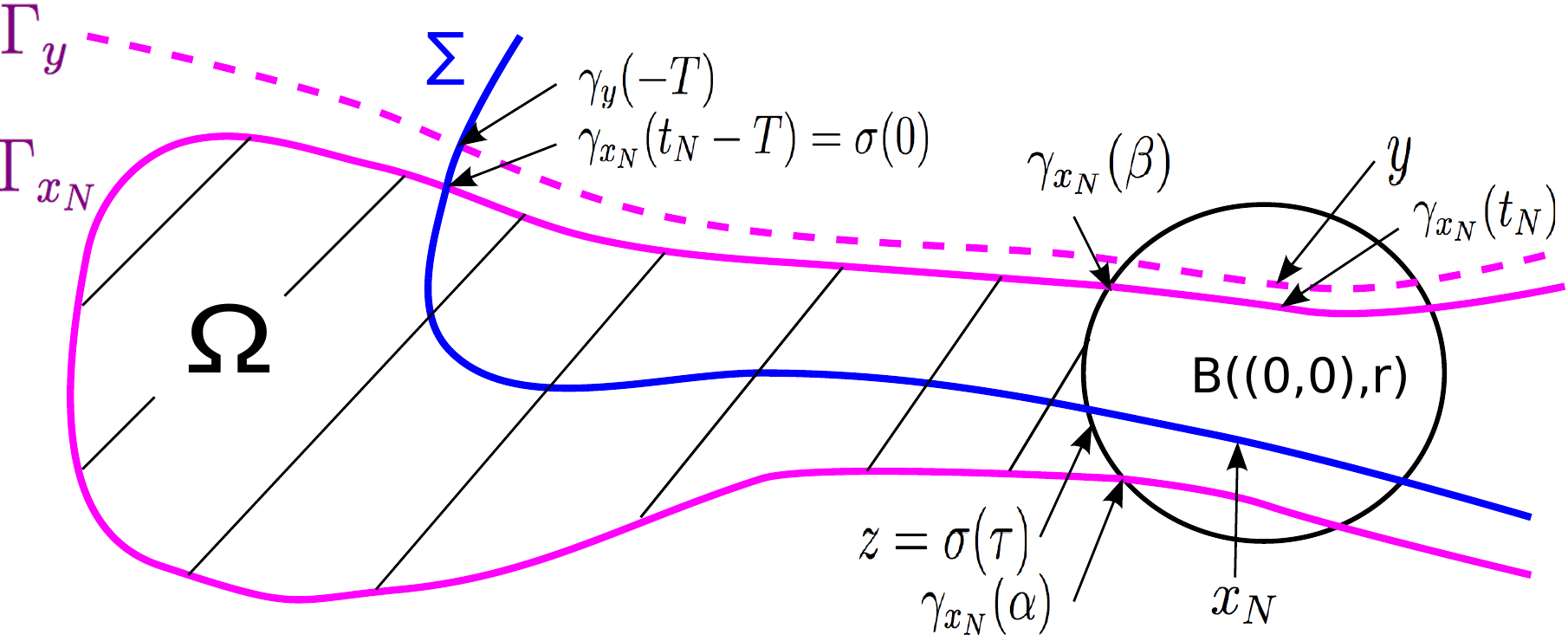}
\caption{The streamlines $\Gamma_{x_N}$ and $\Gamma_y$, the domain $\Omega$, and the trajectory $\Sigma$}
\end{figure}
Let then $\Omega$ be the non-empty bounded domain (open and connected set) surrounded by the closed simple curve $\gamma_{x_N}([\alpha,\beta])\cup A$, where $A$ is the arc on the circle $\partial B(0,r)$ joining $\gamma_{x_N}(\alpha)$ and $\gamma_{x_N}(\beta)$ in such a way that $\Omega\cap B(0,r)=\emptyset$ ($\Omega$ is the hatched region in Figure~1). Let $\Sigma:=\Sigma_{\gamma_{x_N}(t_N-T)}$ be the trajectory of the gradient flow containing the point $\gamma_{x_N}(t_N-T)$, and let $\sigma:=\sigma_{\gamma_{x_N}(t_N-T)}$ be the solution of~\eqref{xix} with initial value $\sigma(0)=\gamma_{x_N}(t_N-T)$. Since $\Sigma$ is orthogonal to $\Gamma_{x_N}$ at the point $\gamma_{x_N}(t_N-T)$ with $t_N-T\in(\alpha,\beta)$ and since $\Sigma$ is unbounded by Lemma~\ref{unboundedbis}, one infers the existence of a real number $\tau\neq 0$ such that
$$\sigma(\tau)\in\partial\Omega\ \hbox{ and }\ \sigma(t)\in\Omega\hbox{ for all }t\in I,$$
where $I$ denotes the open interval $I=(0,\tau)$ if $\tau>0$ (resp. $I=(\tau,0)$ if $\tau<0$).\par
Since $u$ is constant along $\gamma_{x_N}([\alpha,\beta])$ and $t\mapsto u(\sigma(t))$ is increasing by Lemma~\ref{lemsigma}, it follows that $\sigma(\tau)\in A$. Lemma~\ref{lemsigma} also implies that $|u(\sigma(\tau))-u(\sigma(0))|\ge\eta\,|\sigma(\tau)-\sigma(0)|$, that is,
$$|u(\sigma(\tau))-u(\gamma_{x_N}(t_N-T))|\ge\eta\,|\sigma(\tau)-\gamma_{x_N}(t_N-T)|.$$
But $|u(\gamma_{x_N}(t_N-T))|=|u(x_N)|\le M_r=\max_{\overline{B(0,r)}}|u|$ since $|x_N|<r$ by~\eqref{defn0}. Furthermore, $|u(\sigma(\tau))|\le M_r$ too since $\sigma(\tau)\in A\subset\partial B(0,r)$. Thus,
$$\eta\,|\sigma(\tau)-\gamma_{x_N}(t_N-T)|\le|u(\sigma(\tau))-u(\gamma_{x_N}(t_N-T))|\le 2\,M_r.$$
Since $|\sigma(\tau)|=r$ and $|\gamma_{x_N}(t_N-T)|>R$ by~\eqref{defn0}, the triangle inequality yields $\eta\,(R-r)<2\,M_r$. That contradicts the choice of $R$ in~\eqref{rR} and the proof of Lemma~\ref{unboundedter} is thereby complete.\hfill$\Box$\break

The following lemma shows the important property that two streamlines $\Gamma_y$ and $\Gamma_z$ are close to each other in the sense of Hausdorff distance when $y$ and $z$ are close to a given point $x$. We recall that, for any two subsets $A$ and $B$ of $\R^2$, their Hausdorff distance $d_{\mathcal{H}}(A,B)$ is defined as
$$d_{\mathcal{H}}(A,B)=\max\big(\sup\big\{d(a,B);\ a\in A\big\},\sup\big\{d(b,A);\ b\in B\big\}\big),$$
where $d(x,E)=\inf\big\{|x-y|;\ y\in E\big\}$ for any $x\in\R^2$ and $E\subset\R^2$.

\begin{lem}\label{lem2}
For any $x\in\R^2$ and any $\ee>0$, there is $r>0$ such that $d_{\mathcal{H}}(\Gamma_y,\Gamma_z)\le\ee$ for all~$y$ and~$z$ in $B(x,r)$.
\end{lem}

\noindent{\bf{Proof.}} We fix $x\in\R^2$ and $\ee>0$. Up to translation and rotation of the frame, let us assume without loss of generality that~$x$ is  the origin
$$x=0=(0,0),$$
and that $v(x)=v(0)$ points in the $x_2$ direction, that is, $v(0)=|v(0)|\,{\mathrm{e}}_2$ with ${\mathrm{e}}_2=(0,1)$. Firstly, by continuity of $u$, there is $r_0>0$ such that
\beq\label{defr1}
\mathop{\osc}_{B(0,r_0)}u\le\ee\,\eta,
\eeq
where for any non-empty subset $E\subset\R^2$,
$$\mathop{\osc}_Eu=\sup_Eu-\inf_Eu$$
denotes the oscillation of $u$ on $E$, and $\eta>0$ is given in~\eqref{defee0}. Secondly, since
$$\frac{\partial u}{\partial x_1}(0)=v_2(0)=|v(0)|>0$$
and $u$ is (at least) of class $C^1$, there are some real numbers $r\in(0,r_0)$ and~$r^-\le -r<0<r\le r^+$ such that $\frac{\partial u}{\partial x_1}(\cdot,0)>0$ in $(r^-,r^+)$ and
$$\mathop{\bigcup}_{X\in B(0,r)}\Gamma_X=\mathop{\bigcup}_{x_1\in(r^-,r^+)}\Gamma_{(x_1,0)}.$$
In other words, the streamlines of the flow generated by points in the ball $B(0,r)$ are actually generated by the points of the one-dimensional straight curve
$$I=\big\{(x_1,0);\ r^-<x_1<r^+\big\}=(r^-,r^+)\times\{0\},$$
and these streamlines $\Gamma_{(x_1,0)}$ are pairwise distinct since $u$ is one-to-one on $I$.\par
Let now $y$ and $z$ any two points in the ball $B(0,r)$, such that $\Gamma_y\neq\Gamma_z$. From the previous paragraph, there are some real numbers $\tilde{y}_1$ and $\tilde{z}_1$ in $(r^-,r^+)$ such that
$$\Gamma_y=\Gamma_{\tilde{y}}\hbox{ and }\Gamma_z=\Gamma_{\tilde{z}},\ \hbox{ with }\tilde{y}=(\tilde{y}_1,0)\hbox{ and }\tilde{z}=(\tilde{z}_1,0).$$
Since $\Gamma_y\neq\Gamma_z$, there holds $\tilde{y}\neq\tilde{z}$. Without loss of generality, one can then assume that
$$r^-<\tilde{y}_1<\tilde{z}_1<r^+$$
(hence, $u(y)=u(\tilde{y}_1,0)<u(\tilde{z}_1,0)<u(z)$). Denote
\beq\label{defG}
G=\mathop{\bigcup}_{x_1\in(\tilde{y}_1,\tilde{z}_1)}\Gamma_{(x_1,0)}.
\eeq\par
Let us list in this paragraph and the next one some elementary properties of the set $G$. First of all, owing to its definition, $G$ is connected. Secondly, we claim that the set $G$ is open. To show this property, let $X$ be any point in $G$. By definition, there are $x_1\in(\tilde{y}_1,\tilde{z}_1)$ and $\tau\in\R$ such that $X=\gamma_{(x_1,0)}(\tau)$, that is, $(x_1,0)=\gamma_X(-\tau)$. Since $|v(x_1,0)|>0$ and $v_2(x_1,0)=\frac{\partial u}{\partial x_1}(x_1,0)>0$, there is $\rho>0$ such that $\Gamma_{x'}$ intersects the one-dimensional straight curve $(\tilde{y}_1,\tilde{z}_1)\times\{0\}$ for every $x'\in B((x_1,0),\rho)$. On the other hand, by Cauchy-Lipschitz theorem, there is $\delta>0$ such that, for every $x''\in B(X,\delta)$, $\gamma_{x''}(-\tau)\in B((x_1,0),\rho)$ and then $\Gamma_{x''}$ intersects $(\tilde{y}_1,\tilde{z}_1)\times\{0\}$. Therefore, $B(X,\delta)\subset G$ and the set $G$ is then open.\par
Thirdly, we claim that
\beq\label{claimpartial}
\partial G=\Gamma_y\cup\Gamma_z.
\eeq
Indeed, first of all, since the map $X\mapsto\gamma_X(t)$ is continuous for every $t\in\R$ and since $G$ is open, it follows that $\Gamma_y\cup\Gamma_z=\Gamma_{(\tilde{y}_1,0)}\cup\Gamma_{(\tilde{z}_1,0)}\subset\partial G$. Conversely, let now $X$ be any point in $\partial G$. There are then some sequences $(x_{1,n})_{n\in\N}$ in $(\tilde{y}_1,\tilde{z}_1)$ and $(t_n)_{n\in\N}$ in $\R$ such that $\gamma_{(x_{1,n},0)}(t_n)\to X$ as $n\to+\infty$. Since the sequence $((x_{1,n},0))_{n\in\N}$ is bounded, Lemma~\ref{unboundedter} implies that the sequence $(t_n)_{n\in\N}$ is bounded too. Therefore, up to extraction of a subsequence, there holds $x_{1,n}\to x_1\in[\tilde{y}_1,\tilde{z}_1]$ and $t_n\to t\in\R$ as $n\to+\infty$, hence $X=\gamma_{(x_1,0)}(t)$ and $X\in\Gamma_{(x_1,0)}$. But $X\not\in G$ since $X\in\partial G$ and $G$ is open. Thus, either $x_1=\tilde{y}_1$ or $x_1=\tilde{z}_1$. In other words, $X\in\Gamma_{(\tilde{y}_1,0)}\cup\Gamma_{(\tilde{z}_1,0)}=\Gamma_y\cup\Gamma_z$.  Finally, $\partial G\subset\Gamma_y\cup\Gamma_z$ and the claim~\eqref{claimpartial} has been shown.\par
In order to complete the proof of Lemma~\ref{lem2}, consider any point $Y\in\Gamma_y$. Let $\Sigma'_Y$ be the restriction to $\overline{G}$ of the trajectory $\Sigma_Y$ of the gradient flow $\dot{\sigma}=\nabla u(\sigma)$ containing $Y\,(\in\Gamma_y\subset\partial G)$ and let $\sigma_Y$ be the solution of~\eqref{xix} starting at $\sigma_Y(0)=Y$. Remember that the function $t\mapsto g_Y(t):=u(\sigma_Y(t))$ is increasing in $\R$ by Lemma~\ref{lemsigma}. Since $u(Y)=u(y)=u(\tilde{y}_1,0)<u(\tilde{z}_1,0)=u(z)$ and $\nabla u(Y)$ is orthogonal to~$\partial G$ at $Y$, there is $t_0>0$ such that $\sigma_Y(t)\in G$ for all $t\in(0,t_0)$. Furthermore, one infers from the definition of the streamlines and from~\eqref{defG} that $u$ is bounded in the open set $G$. Since $\partial G=\Gamma_y\cup\Gamma_z$ and $g'_Y(t)\ge\eta^2>0$ in $\R$, there is $t_1>0$ such that $\sigma_Y(t)\in G$ for all $t\in(0,t_1)$ and $Z:=\sigma_Y(t_1)\in\partial G$ (hence, $Z\in\Gamma_z$ since $u(\sigma_Y(\cdot))$ is increasing). Therefore, $\Sigma'_Y$ can be parametrized by $\sigma_Y(t)$ for $t\in[0,t_1]$ and, by Lemma~\ref{lemsigma},
$$|Y-Z|=|\sigma_Y(0)-\sigma_Y(t_1)|\le\frac{u(\sigma_Y(t_1))-u(\sigma_Y(0))}{\eta}=\frac{u(Z)-u(Y)}{\eta}=\frac{u(z)-u(y)}{\eta}.$$
Since both points $y$ and $z$ belong to $B(0,r)\subset B(0,r_0)$ and $\mathop{\osc}_{B(0,r_0)}u\le\ee\,\eta$ by~\eqref{defr1}, one infers that $|Y-Z|\le\ee$. Since $Z\in\Gamma_z$, it follows that
$$d(Y,\Gamma_z)\le\ee.$$
As~$Y$ was arbitrary in $\Gamma_y$ and both points $y$ and $z$ play a similar role, one concludes that $d_{\mathcal{H}}(\Gamma_y,\Gamma_z)\le\ee$. The proof of Lemma~\ref{lem2} is thereby complete.\hfill$\Box$\break

Based on the previous results, we show in the following lemma that the level curves of $u$ foliate the plane $\R^2$ in a monotone way.

\begin{lem}\label{foliation}
For any trajectory $\Sigma$ of the gradient flow, there holds
\beq\label{SigmaGamma}
\bigcup_{x\in\Sigma}\Gamma_x=\R^2.
\eeq
\end{lem}

Notice that this lemma implies that any level set of $u$ has only one connected component. Indeed, for any $\lambda\in\R$, Lemma~\ref{unboundedbis} implies that there is a unique $x\in\Sigma$ such that $u(x)=\lambda$. Since $u$ is constant along any streamline, it follows then from Lemma~\ref{foliation} that the level set $\big\{y\in\R^2;\ u(y)=\lambda\big\}$ is equal to the streamline~$\Gamma_x$ and thus has only one connected component.\hfill\break

\noindent{\bf{Proof.}} Consider any trajectory $\Sigma$ of the gradient flow $\dot\sigma=\nabla u(\sigma)$ and let $X\in\Sigma$, that is, $\Sigma=\Sigma_X$. Denote
$$E=\bigcup_{x\in\Sigma}\Gamma_x$$
and let us show that $E=\R^2$. First of all, $E$ is not empty and, with the same arguments as for the set $G$ defined in~\eqref{defG}, one gets that $E$ is open.\par
To conclude that $E=\R^2$, it is then sufficient to show that $E$ is closed. So, let $(y_n)_{n\in\N}$ be a sequence in $E$ and $y\in\R^2$ such that $y_n\to y$ as $n\to+\infty$. Owing to the definition of $E$ and given the parametrizations $t\mapsto\gamma_x(t)$ of $\Gamma_x$ for every $x\in\R^2$, there are some sequences $(x_n)_{n\in\N}$ in $\Sigma$ and $(t_n)_{n\in\N}$ in $\R$ such that
$$y_n=\gamma_{x_n}(t_n)\ \hbox{ for all }n\in\N.$$
Furthermore, since $\sigma_X(\cdot)$ is a parametrization of $\Sigma$, there is a sequence $(\tau_n)_{n\in\N}$ in $\R$ such that
$$x_n=\sigma_X(\tau_n)\ \hbox{ for all }n\in\N.$$
Since $u(\sigma_X(\tau_n))=u(x_n)=u(y_n)\to u(y)$ as $n\to+\infty$ (by continuity of $u$ and definition of the streamlines $\Gamma_{x_n}$) and $|u(\sigma_X(\tau))|\to+\infty$ as $|\tau|\to+\infty$ by Lemma~\ref{unboundedbis}, it follows that the sequence $(\tau_n)_{n\in\N}$ is bounded. Up to extraction of a subsequence, there is $\tau\in\R$ such that $\tau_n\to\tau$ as $n\to+\infty$, hence
$$x_n=\sigma_X(\tau_n)\to x:=\sigma_X(\tau)\,(\in\Sigma)\ \hbox{ as }n\to+\infty.$$
Consider now any $\ee>0$. For $n$ large enough, one has $|x_n-x|\le\ee$. Moreover, from Lemma~\ref{lem2}, there holds
$$d_{\mathcal{H}}(\Gamma_{y_n},\Gamma_y)\le\ee\ \hbox{ for }n\hbox{ large enough},$$
that is, $d_{\mathcal{H}}(\Gamma_{x_n},\Gamma_y)\le\ee$ ($\Gamma_{x_n}=\Gamma_{y_n}$ by definition of $x_n$). Therefore, $d(x_n,\Gamma_y)\le\ee$ for $n$ large enough, hence $d(x,\Gamma_y)\le2\ee$. As a consequence, since $\ee>0$ can be arbitrarily small, $d(x,\Gamma_y)=0$. On the other hand, $\Gamma_y$ is a closed subset of $\R^2$ from its definition and from Lemma~\ref{unbounded}. Thus, $x\in\Gamma_y$. In other words, $\Gamma_y=\Gamma_x$ and $y\in\Gamma_x$. Finally, $y\in E$ and $E$ is closed. As a conclusion, $E=\R^2$ and the proof of Lemma~\ref{foliation} is thereby complete.\hfill$\Box$

\begin{rem}{\rm As a immediate corollary of Lemma~\ref{foliation}, it follows that the trajectories of the gradient flow foliate the whole plane $\R^2$ in the sense that the family of trajectories of the gradient flow can be parametrized by the points along any streamline of the flow. More precisely, for any streamline $\Gamma$ of the flow, there holds
$$\bigcup_{y\in\Gamma}\Sigma_y=\R^2.$$
The property will actually not be used in the sequel, but we state it as an interesting counterpart of~\eqref{SigmaGamma}.}
\end{rem}

From Lemma~\ref{foliation}, the level curves of $u$ foliate the plane $\R^2$ in the sense that the family of streamlines can be parametrized by the points along any trajectory of the gradient flow. Since $\Delta u$ turns out to be constant along any streamline, the function $u$ will then satisfy a simple semilinear elliptic equation in $\R^2$. Namely, the following result holds.

\begin{lem}\label{elliptic}
There is a $C^1$ function $f:\R\to\R$ such that $u$ is a classical solution of
\beq\label{eqelliptic}
\Delta u+f(u)=0\ \hbox{ in }\R^2.
\eeq
\end{lem}

\noindent{\bf{Proof.}} Let $\Sigma=\Sigma_0$ be the trajectory of the gradient flow going through the origin, and let $\sigma=\sigma_0:\R\to\R^2$ be its $C^1(\R)$ parametrization defined by~\eqref{xix} with $x=0=(0,0)$. As already underlined in the proof of Lemmas~\ref{unboundedbis} and~\ref{lemsigma}, the $C^1(\R)$ function $g:\R\to\R,\ t\mapsto g(t):=u(\sigma(t))$ is increasing and $g'(t)\ge\eta^2>0$ for all $t\in\R$. Let $g^{-1}\in C^1(\R)$ be the reciprocal function of $g$ and let us now define
$$\baa{rcl}
f:\R & \!\!\to\!\! & \R\vspace{3pt}\\
s & \!\!\mapsto\!\! & f(s):=-\Delta u\big(\sigma(g^{-1}(s))\big),\eaa$$
that is,
$$f(g(t))=-\Delta u(\sigma(t))$$
for all $t\in\R$. Since $u$ is of class $C^3(\R^2)$ and both $\sigma$ and $g^{-1}$ are of class $C^1(\R)$, one infers that $f$ is of class $C^1(\R)$ too.\par
Let us then show that $u$ is a classical solution of the elliptic equation~\eqref{eqelliptic}. Indeed, since
$$\Delta u=\frac{\partial v_2}{\partial x_1}-\frac{\partial v_1}{\partial x_2}\ \hbox{ in }\R^2$$
and since this $C^1(\R^2)$ function satisfies
$$v\cdot\nabla(\Delta u)=v\cdot\nabla\Big(\frac{\partial v_2}{\partial x_1}-\frac{\partial v_1}{\partial x_2}\Big)=0\ \hbox{ in }\R^2$$
by~\eqref{1}, one infers that the function $\Delta u$ is constant along any streamline of $v$, that is, along any level curve of $u$.\par
Let finally $x$ be any point in $\R^2$. From Lemma~\ref{foliation}, the streamline $\Gamma_x$ intersects $\Sigma$. Therefore, since $u$ is constant along $\Gamma_x$ and $\sigma:\R\to\Sigma$ is one-to-one (as $g=u\circ\sigma:\R\to\R$ is one-to-one too), there is a unique $t_x\in\R$ such that $\sigma(t_x)\in\Gamma_x$, and
$$g(t_x)=u(\sigma(t_x))=u(x).$$
As a conclusion, since the function $\Delta u$ is constant on the streamline $\Gamma_x$ (containing both $x$ and $\sigma(t_x)$), one infers from the definitions of $g$ and $f$ that
$$
\Delta u(x)=\Delta u(\sigma(t_x))=-f(g(t_x))=-f(u(x)).$$
The proof of Lemma~\ref{elliptic} is thereby complete.\hfill$\Box$

\begin{rem}{\rm If $v$ is not assumed to be in $L^{\infty}(\R^2)$ anymore, then $v$ is still locally bounded (since it is at least continuous). In that case, for every $x\in\R^2$, the function $\sigma_x$ solving~\eqref{xix} would be defined in a maximal interval $(t^-_x,t^+_x)$ with $-\infty\le t^-_x<0<t^+_x\le+\infty$, and $|\sigma_x(t)|\to+\infty$ as $t\to t^-_x$ if $t^-_x\in\R$ (resp. as $t\to t^+_x$ if $t^+_x\in\R$). Furthermore, the arguments used in the proof of Lemma~\ref{unboundedbis} still imply that $|\sigma_x(t)|\to+\infty$ as $t\to t^-_x$ if $t^-_x=-\infty$ (resp. as $t\to t^+_x$ if $t^+_x=+\infty$). In particular, any trajectory $\Sigma$ of the gradient flow is still unbounded. Property~\eqref{ualphabeta} in Lemma~\ref{lemsigma} still holds as well, as soon as $I=\big[\min(\alpha,\beta),\max(\alpha,\beta)\big]\subset(t^-_x,t^+_x)$, for any $x\in\Sigma$ with $\sigma=\sigma_x$. For every $x\in\R^2$, the arc length parametrization $t\mapsto\varsigma_x(t)$ of $\Sigma_x$ solving
$$\dot\varsigma_x(t)=\frac{\nabla u(\varsigma_x(t))}{|\nabla u(\varsigma_x(t))|}$$
is defined in the whole interval $\R$, and the function $h:t\mapsto u(\varsigma_x(t))$ satisfies $h'(t)=|\nabla u(\varsigma_x(t))|\ge\eta>0$, hence $u(\varsigma_x(t))\to\pm\infty$ as $t\to\pm\infty$.\par
Similarly, still if $v$ is not assumed to be in $L^{\infty}(\R^2)$ anymore, for every $x\in\R^2$, the function $\gamma_x$ solving~\eqref{defgammax} would be defined in a maximal interval $(\tau^-_x,\tau^+_x)$ with $-\infty\le\tau^-_x<0<\tau^+_x\le+\infty$, and $|\gamma_x(t)|\to+\infty$ as $t\to\tau^-_x$ if $\tau^-_x\in\R$ (resp. as $t\to\tau^+_x$ if $\tau^+_x\in\R$. Moreover, if $\tau^-_x=-\infty$ (resp. $\tau^+_x=+\infty$), the arguments used in the proof of Lemma~\ref{unbounded} still imply that $|\gamma_x(t)|\to+\infty$ as $t\to\tau^-_x=-\infty$ (resp. as $t\to\tau^+_x=+\infty$). Therefore, in all cases, whether $\tau^{\pm}_x$ be finite or not, one has $|\gamma_x(t)|\to+\infty$ as $t\to\tau^{\pm}_x$ and the streamline $\Gamma_x$ is unbounded. In particular, for every $x\in\R^2$, the arc length parametrization $t\mapsto\zeta_x(t)$ of $\Gamma_x$ solving
$$\dot\zeta_x(t)=\frac{v(\zeta_x(t))}{|v(\zeta_x(t))|}$$
is defined in the whole interval $\R$. With the unboundedness of each trajectory of the gradient flow and with~\eqref{ualphabeta}, the arguments used in the proof of Lemma~\ref{unboundedter} can be repeated with the parametrizations $\zeta_x$ instead of $\gamma_x$: in other words, for every bounded set $K$, there holds $|\zeta_x(t)|\to+\infty$ as $|t|\to+\infty$ uniformly in $x\in K$. Similarly, the proofs of Lemmas~\ref{lem2},~\ref{foliation} and~\ref{elliptic} can be done with the parametrizations $\zeta_x$ and $\varsigma_x$ instead of $\gamma_x$ and $\sigma_x$.\par
To sum up, the conclusions of the aforementioned lemmas still hold if $v$ is not assumed to be in $L^{\infty}(\R^2)$. Actually, the purpose of this remark is to make a connection with the beginning of the proof of~\cite[Theorem~1.1]{hn}, where the flow $v$, which was there defined in a two-dimensional strip, was indeed not assumed to be a priori bounded. Namely, the beginnings of the proofs of~\cite[Theorem~1.1]{hn} and of Theorem~\ref{th1} of the present paper are similar, even if more details and additional properties are proved here, such as the unboundedness of the trajectories of the gradient flow and the uniform unboundedness of the streamlines emanating from a bounded region. However, the remaining part of the proof of Theorem~\ref{th1} of the present paper, as well as the proof of Proposition~\ref{growth} below, strongly uses the boundedness of $v$ (and, as emphasized in Section~\ref{intro}, the conclusion of Theorem~\ref{th1} is not valid in general without the boundedness of $v$).}
\end{rem}


\subsection{End of the proof of Theorem~\ref{th1} in the general case}\label{sec22}

In order to complete the proof of Theorem~\ref{th1}, the following propositions provide key estimates on the oscillations of the argument of the vector field $v$. These estimates, which will be used for scaled or shifted fields, are thus established for general solutions $w$ of the Euler equations~\eqref{1}. To state these estimates, let us first introduce a few more notations.\par
For any $C^2(\R^2)$ solution $w$ of~\eqref{1} (associated with a pressure $p=p_w$) and satisfying~\eqref{defee0} for some $\eta\in(0,1]$, that is,
\beq\label{etaw}
0<\eta\le|w(x)|\le\eta^{-1}\ \hbox{ for all }x\in\R^2,
\eeq
there is a $C^2(\R^2)$ function $\phi_w$ such that
\beq\label{argv}
\frac{w(x)}{|w(x)|}=(\cos\phi_w(x),\sin\phi_w(x))\ \hbox{ for all }x\in\R^2.
\eeq
This function, which is the argument of the flow $w$, is uniquely defined in $\R^2$ up to an additive constant which is a multiple of $2\pi$. Its oscillation
$$\mathop{\osc}_E\phi_w=\sup_E\phi_w-\inf_E\phi_w$$
in any non-empty subset $E\subset\R^2$ is uniquely defined (namely it does not depend on the choice of this additive constant). Similarly, after denoting $u_w$ the unique stream function of the flow $w$ (uniquely defined by $\nabla^\perp u_w=w$ and $u_w(0)=0$), there is a $C^2(\R^2)$ function $\varphi_w:\R^2\to\R$ such that
\beq\label{defvarphi}
\frac{\nabla u_w(x)}{|\nabla u_w(x)|}=(\cos\varphi_w(x),\sin\varphi_w(x))
\eeq
for all $x\in\R^2$. Lastly, there is an integer $q\in\Z$ such that
\beq\label{phivarphi}
\varphi_w(x)=\phi_w(x)-\pi/2+2\pi q\ \hbox{ for all }x\in\R^2.
\eeq
In particular,
\beq\label{oscvarphi}
\mathop{\osc}_E\varphi_w=\mathop{\osc}_E\phi_w
\eeq
for every non-empty subset $E\subset\R^2$.

The key-estimate is the following logarithmic upper bound of the oscillations of the arguments of the solutions $w$ of~\eqref{1} and~\eqref{etaw} in large balls, given an upper bound in smaller balls.

\begin{pro}\label{growth}
For any $\eta\in(0,1]$, there is a positive real number $C_\eta$ such that, for any $C^2(\R^2)$ solution $w$ of~\eqref{1} satisfying~\eqref{etaw} and for any $R\ge2$, if
$$\mathop{\osc}_{B(x,1)}\phi_w<\frac{\pi}{4}\ \hbox{ for all }x\in B(0,R),$$
then
\beq\label{loggrowth}
\mathop{\osc}_{B(0,R)}\phi_w\le C_\eta\ln R.
\eeq
\end{pro}

In order not to loose the main thread of the proof of Theorem~\ref{th1}, the proof of Proposition~\ref{growth} is postponed in Section~\ref{sec3}.\par
The second key-estimate in the proof of Theorem~\ref{th1} is the following lower bound of the oscillations of the arguments of the solutions $w$ of~\eqref{1} and~\eqref{etaw} in some balls of radius $1/2$, given a lower bound in the unit ball.

\begin{pro}\label{pro2}
For any $\eta\in(0,1]$, there is a positive real number $R_\eta$ such that, for any $C^2(\R^2)$ solution $w$ of~\eqref{1} satisfying~\eqref{etaw}, if
$$\mathop{\osc}_{B(0,1)}\phi_w\ge\frac{\pi}{4},$$
then there is a point $x\in B(0,R_\eta)$ such that
$$\mathop{\osc}_{B(x,1/2)}\phi_w\ge\frac{\pi}{4}.$$
\end{pro}

\noindent{\bf{Proof.}} It is based on Proposition~\ref{growth}, on some scaling arguments, on the derivation of some uniformly elliptic linear equations satisfied by the arguments $\phi_w$ of the $C^2(\R^2)$ solutions $w$ of~\eqref{1} satisfying~\eqref{etaw}, and on some results of Moser~\cite{m} on the solutions of such elliptic equations.\par
Consider any $C^2(\R^2)$ solution $w$ of~\eqref{1} and~\eqref{etaw}, and let us first derive a linear elliptic equation for its argument $\phi_w$ (see also~\cite{f1,kn} for the derivation of such equations). Let $x$ be any point in $\R^2$. Assume first that $\nabla u_w(x)$ is not parallel to the vector $(0,1)$, that is $\frac{\partial u_w}{\partial x_1}(x)\neq 0$. In other words, by continuity of $\nabla u_w$, there is $k\in\N$ such that
$$\varphi_w=\arctan\Big(\frac{\partial u_w}{\partial x_2}/\frac{\partial u_w}{\partial x_1}\Big)+k\pi$$
in a neighborhood of $x$. Hence, a straightforward calculation leads to
\beq\label{varphi}
|\nabla u_w|^2\,\nabla\varphi_w=\frac{\partial u_w}{\partial x_1}\,\nabla\Big(\frac{\partial u_w}{\partial x_2}\Big)-\frac{\partial u_w}{\partial x_2}\,\nabla\Big(\frac{\partial u_w}{\partial x_1}\Big)
\eeq
in a neighborhood of $x$. Similarly, if $\nabla u_w(x)$ is not parallel to the vector $(1,0)$, then $\cot(\varphi_w)=\frac{\partial u_w}{\partial x_1}/\frac{\partial u_w}{\partial x_2}$ in a neighborhood of $x$, and formula~\eqref{varphi} still holds in a neighborhood of $x$. Therefore,~\eqref{varphi} holds in $\R^2$. On the other hand, by Lemma~\ref{elliptic} and by differentiating with respect to both variables $x_1$ and~$x_2$ the elliptic equation~\eqref{eqelliptic} satisfied by~$u_w$ (formula~\eqref{eqelliptic} holds for some $C^1(\R)$ function $f=f_w$ depending on $w$), it follows that
\beq\label{equx}
\Delta\Big(\frac{\partial u_w}{\partial x_1}\Big)+f_w'(u_w)\,\frac{\partial u_w}{\partial x_1}=\Delta\Big(\frac{\partial u_w}{\partial x_2}\Big)+f_w'(u_w)\,\frac{\partial u_w}{\partial x_2}=0\ \hbox{ in }\R^2.
\eeq
Together with~\eqref{varphi}, one obtains that ${\div}\big(|\nabla u_w|^2\nabla\varphi_w\big)=0$ in $\R^2$. In other words, thanks to~\eqref{phivarphi} and $\nabla^\perp u_w=w$, there holds
\beq\label{eqargv}
{\div}\big(|w|^2\,\nabla\phi_w\big)=0\ \hbox{ in }\R^2.
\eeq\par
Since $w$ satisfies the uniform bounds~\eqref{etaw}, it then follows from the proof of~\cite[Theorem~4]{m} that there are some positive real numbers $M_\eta$ and $\alpha_\eta$, depending on $\eta$ only (and not on $w$ and $\phi_w$) such that
$$\mathop{\osc}_{B(0,R)}\phi_w\ge M_\eta\,R^{\alpha_\eta}\!\mathop{\osc}_{B(0,1)}\phi_w\ \hbox{ for all }R\ge1.$$
There exists then a positive real number $R_\eta$, depending on $\eta$ only, such that
$$R_\eta\ge1\ \hbox{ and }\ M_\eta R_\eta^{\alpha_\eta}>\frac{4C_\eta}{\pi}\ln(2R_\eta),$$
where $C_\eta$ is the positive constant given in Proposition~\ref{growth}. The previous two formulas yield
\beq\label{Retabis}
\mathop{\osc}_{B(0,R_\eta)}\phi_w>\frac{4C_\eta}{\pi}\,\ln(2R_\eta)\,\mathop{\osc}_{B(0,1)}\phi_w.
\eeq\par
Let us now check that Proposition~\ref{pro2} holds with this value $R_\eta$. Consider any $C^2(\R^2)$ solution $w$ of~\eqref{1} (with pressure $p_w$) satisfying~\eqref{etaw} and ${\osc}_{B(0,1)}\phi_w\ge\pi/4$. Assume by way of contradiction that the conclusion does not hold, that is,
\beq\label{oscphiw}
\mathop{\osc}_{B(x,1/2)}\phi_w<\frac{\pi}{4}\ \hbox{ for all }x\in B(0,R_\eta).
\eeq
Define now
$$\tilde{w}(x)=w\Big(\frac{x}{2}\Big)\ \hbox{ for }x\in\R^2.$$
The $C^2(\R^2)$ function $\tilde{w}$ still satisfies~\eqref{etaw}, as well as~\eqref{1} with pressure $\tilde{p}(x)=p_w(x/2)$. The arguments $\phi_{\tilde{w}}$ and $\phi_w$ of $\tilde{w}$ and $w$ satisfy $\phi_{\tilde{w}}(x)=\phi_w(x/2)+2k\pi$ in $\R^2$, for some integer $k\in\Z$, hence property~\eqref{oscphiw} translates into
$$\mathop{\osc}_{B(x,1)}\phi_{\tilde{w}}<\frac{\pi}{4}\ \hbox{ for all }x\in B(0,2R_\eta).$$
Proposition~\ref{growth} applied with $\tilde{w}$ and $R:=2R_\eta\ge2$ then yields
$$\mathop{\osc}_{B(0,2R_\eta)}\phi_{\tilde{w}}\le C_\eta\ln(2R_\eta),$$
that is, ${\osc}_{B(0,R_\eta)}\phi_w\le C_\eta\ln(2R_\eta)$. On the other hand, it follows from~\eqref{Retabis} and the assumption ${\osc}_{B(0,1)}\phi_w\ge\pi/4$ that
$$\mathop{\osc}_{B(0,R_\eta)}\phi_w>C_\eta\,\ln(2R_\eta).$$
One has then reached a contradiction, and the proof of Proposition~\ref{pro2} is thereby complete.\hfill$\Box$

\begin{rem}{\rm From the derivation of~\eqref{eqargv} in the proof of Proposition~\ref{pro2}, we point out that, for any $C^2(\Omega)$ flow $w$ solving~\eqref{1} in a domain $\Omega\subset\R^2$, the equation
$${\div}\big(|w|^2\,\nabla\phi_w\big)=0$$
holds in a neighborhood of any point $x$ where $|w(x)|\neq0$, where $\phi_w$ is a $C^2$ argument of the flow $w$ given by~\eqref{argv} in a neighborhood of $x$. That equation also holds globally in $\Omega$ if $\Omega$ is simply connected and $w$ has no stagnation point in the whole domain~$\Omega$.}
\end{rem}

With Proposition~\ref{pro2} and its proof in hand, the proof of Theorem~\ref{th1} can then be carried out.\hfill\break

\noindent{\bf{Proof of Theorem~\ref{th1}.}} Let $v$ be a $C^2(\R^2)$ solution of~\eqref{1} satisfying~\eqref{defee0} with $0<\eta\le1$, and associated with a pressure $p$. Let $\phi=\phi_v$ be the $C^2(\R^2)$ argument of $v$, satisfying~\eqref{argv} with $v$ instead of $w$. We want to show that $\phi$ is constant. Assume by way of contradiction that $\phi$ is not constant. Then, since $\phi$ satisfies the elliptic equation ${\div}(|v|^2\nabla\phi)=0$ in $\R^2$ with $0<\eta^2\le|v|^2\le\eta^{-2}$ in $\R^2$, it follows in particular from~\cite[Theorem~4]{m} that ${\osc}_{B(0,\rho)}\phi\to+\infty$ (and even grows algebraically) as $\rho\to+\infty$. Therefore, there is $\rho>0$ such that
$$\mathop{\osc}_{B(0,\rho)}\phi\ge\frac{\pi}{4}.$$
Define now
$$w_1(x)=v(\rho x)\ \hbox{ for }x\in\R^2.$$
The $C^2(\R^2)$ vector field $w_1$ still solves~\eqref{1} (with pressure $p_{w_1}(x)=p(\rho x)$) and satisfies~\eqref{etaw}. Furthermore, its argument $\phi_{w_1}$ satisfies ${\osc}_{B(0,1)}\phi_{w_1}\ge\pi/4$, because $\phi_{w_1}(x)=\phi(\rho x)+2k_1\pi$ in $\R^2$ for some $k_1\in\Z$. It then follows from Proposition~\ref{pro2} that there is a point $y_1\in B(0,R_\eta)$ such that
$$\mathop{\osc}_{B(y_1,1/2)}\phi_{w_1}\ge\frac{\pi}{4},$$
where $R_\eta>0$ is given by Proposition~\ref{pro2} and depends on $\eta$ only. The above formula means that
$$\mathop{\osc}_{B(\rho y_1,\rho/2)}\phi\ge\frac{\pi}{4}.$$\par
Define now
$$w_2(x)=w_1\Big(y_1+\frac{x}{2}\Big).$$
The $C^2(\R^2)$ vector field $w_2$ still solves~\eqref{1} (with pressure $p_{w_2}(x)=p_{w_1}(y_1+x/2)$) and satisfies~\eqref{etaw}. Furthermore, its argument $\phi_{w_2}$ satisfies ${\osc}_{B(0,1)}\phi_{w_2}\ge\pi/4$, because $\phi_{w_2}(x)=\phi_{w_1}(y_1+x/2)+2k_2\pi$ in $\R^2$ for some $k_2\in\Z$. It then follows from Proposition~\ref{pro2} that there is a point $y_2\in B(0,R_\eta)$ such that ${\osc}_{B(y_2,1/2)}\phi_{w_2}\ge\pi/4$. This means that ${\osc}_{B(y_1+y_2/2,1/4)}\phi_{w_1}\ge\pi/4$, that is,
$$\mathop{\osc}_{B(\rho y_1+\rho y_2/2,\rho/4)}\phi\ge\frac{\pi}{4}.$$
By an immediate induction, there exists a sequence of points $(y_n)_{n\in\N}$ of the ball $B(0,R_\eta)$ such that
\beq\label{oscBn}
\mathop{\osc}_{B(\rho y_1+\rho y_2/2+\cdots+\rho y_n/2^{n-1},\rho/2^n)}\phi\ge\frac{\pi}{4}
\eeq
for every $n\in\N$. There is then a point $\overline{z}\in\R^2$ such that $z_n:=\rho y_1+\rho y_2/2+\cdots+\rho y_n/2^{n-1}\to\overline{z}$ as $n\to+\infty$. But since $\phi$ is (at least) continuous at $\overline{z}$, there is $\epsilon>0$ such that ${\osc}_{B(\overline{z},\epsilon)}\phi<\pi/4$, hence ${\osc}_{B(z_n,\rho/2^n)}\phi\le{\osc}_{B(\overline{z},\epsilon)}<\pi/4$ for all $n$ large enough so that $B(z_n,\rho/2^n)\subset B(\overline{z},\epsilon)$. This contradicts~\eqref{oscBn}.\par
As a conclusion, the argument $\phi$ of $v$ is constant, which means that $v(x)$ is parallel to the constant vector $e=(\cos\phi,\sin\phi)$ for all $x\in\R^2$. In other words, $v$ is a shear flow and, since it is divergence free, it can then be written as in~\eqref{shear}, namely $v(x)=V(x\cdot e^{\perp})\,e$. Lastly, by continuity and~\eqref{hyp1}, the function $V:\R\to\R$ has a constant strict sign. The proof of Theorem~\ref{th1} is thereby complete.\hfill$\Box$


\subsection{End of the proof of Theorem~\ref{th1} in the monotone case}\label{sec23}

The goal of this section is to provide an alternate proof of Theorem~\ref{th1}, without making use of Proposition~\ref{growth}, but with the additional assumption that $v(x)$ belongs to a fixed half-plane $\big\{x\in\R^2;\ x\cdot e>0\big\}$ for some unit vector $e\in\mathbb{S}^1$, see Remark~\ref{remmonotone}. That is, we assume in this section that
\beq\label{hypmonotone}
v\cdot e>0\ \hbox{ in }\R^2.
\eeq\par
Of course, the assumption~\eqref{hypmonotone} implies that the $C^2(\R^2)$ argument $\phi$ of $v$ defined by~\eqref{argv} (with $v$ instead of $w$) is bounded. As it solves the elliptic equation ${\div}(|v|^2\nabla\phi)=0$ in $\R^2$ with $0<\eta^2\le|v|^2\le\eta^{-2}$ in $\R^2$, it follows immediately from~\cite[Theorem~4]{m} that $\phi$ is constant, hence $v$ is a shear flow.\par
But here, with the assumption~\eqref{hypmonotone}, we would like to provide an alternate proof of the main result, without making use of the argument $\phi$ of $v$. To do so, up to rotation of the frame, let us assume without loss of generality that $e=(-1,0)$. In other words, we assume here that
$$\frac{\partial u}{\partial x_2}>0\ \hbox{ in }\R^2.$$
From Lemma~\ref{elliptic}, both $C^2(\R^2)$ functions $\frac{\partial u}{\partial x_1}$ and $\frac{\partial u}{\partial x_2}$ satisfy the equation~\eqref{equx} (here, $u_w=u$ and $f_w=f$), with $f'(u(\cdot))\in L^{\infty}_{loc}(\R^2)$. Since $|\nabla u|=|v|$ is bounded and $\frac{\partial u}{\partial x_2}$ is positive, it follows directly from~\cite[Theorem~1.8]{bcn} that $\frac{\partial u}{\partial x_1}$ and $\frac{\partial u}{\partial x_2}$ are proportional. Hence, $\nabla u$ is parallel to a constant vector, that is, $v$ is parallel to a constant vector, and the conclusion follows as in Section~\ref{sec22} above.


\SE{Proof of Proposition~\ref{growth}}\label{sec3}

As already emphasized in Section~\ref{sec22}, the proof of Theorem~\ref{th1} relies on Proposition~\ref{growth}. The main underlying idea in the proof of Proposition~\ref{growth} is the fact that, if the arguments $\phi_w$ and $\varphi_w$ of $w$ and $\nabla u_w$ given in~\eqref{argv} and~\eqref{defvarphi} oscillate too much in some large balls, then there would be some trajectories $\Sigma$ of the gradient flow $\dot\sigma=\nabla u_w(\sigma)$ which would turn too many times. Since $u_w$ grows at least (and at most too) linearly along the trajectories $\Sigma$ while it grows at most linearly in any direction, one would then get a contradiction.\par
The proof of Proposition~\ref{growth} is itself divided into several lemmas and subsections. In Section~\ref{sec31}, we show some estimates on the oscillation of the argument of any $C^1$ embedding between two points, in terms of its surrounding arcs. In Sections~\ref{sec32} and~\ref{sec33}, we show the logarithmic growth of the argument $\varphi_w$ of $\nabla u_w$ along the trajectories of the gradient flow and along the streamlines of the flow. Lastly, we complete the proof of Proposition~\ref{growth} in Section~\ref{sec34}.


\subsection{Oscillations of the argument of any $C^1$ embedding}\label{sec31}

Let us first introduce a few useful notations. For any two points $x\neq y$ in $\R^2$, let
$$L_{x,y}=\big\{(1-t)x+ty;\ t\in\R\big\}$$
be the line containing both $x$ and $y$. Let
$$[x,y]=\big\{(1-t)x+ty;\ t\in[0,1]\big\}$$
be the segment joining $x$ and $y$. Similarly, we denote $(x,y]=\big\{(1-t)x+ty;\ t\in(0,1]\big\}$, $[x,y)=\big\{(1-t)x+ty;\ t\in[0,1)\big\}$ and $(x,y)=\big\{(1-t)x+ty;\ t\in(0,1)\big\}$. We also say that a point $z\in L_{x,y}$ is on the left (resp. right) of $x$ with respect to $y$ if $(z-x)\cdot(y-x)<0$ (resp. $>0$).\par
For any $C^1$ embedding $\xi:\R\to\R^2$ (one-to-one map such that $|\dot\xi(t)|>0$ for all $t\in\R$), let $\theta$ be one of its continuous arguments, defined by
\beq\label{deftheta}
\frac{\dot\xi(t)}{|\dot\xi(t)|}=\big(\!\cos\theta(t),\sin\theta(t)\big)
\eeq
for all $t\in\R$. The continuous function $\theta$ is well defined, and it is unique up to multiples of $2\pi$ (it is unique if it is normalized so that $\theta(0)\in[0,2\pi)$). The following two lemmas provide some fundamental estimates of the oscillations of the argument $\theta(t)$ between any two real numbers $a$ and $b$, according to the number of times $\xi(t)$ turns around $\xi(a)$ and $\xi(b)$ for $t\in[a,b]$.\par
First of all, for any $a<b\in\R$, we say that the arc $\xi([a,b])$ is non-intersecting if
\beq\label{defnoninter}
\xi([a,b])\cap(\xi(a),\xi(b))=\emptyset\ \hbox{ or }\ \xi([a,b])\cap\big(L_{\xi(a),\xi(b)}\backslash[\xi(a),\xi(b)]\big)=\emptyset,
\eeq
see Figure~2.
\begin{figure}\centering
\subfigure{\includegraphics[scale=0.9]{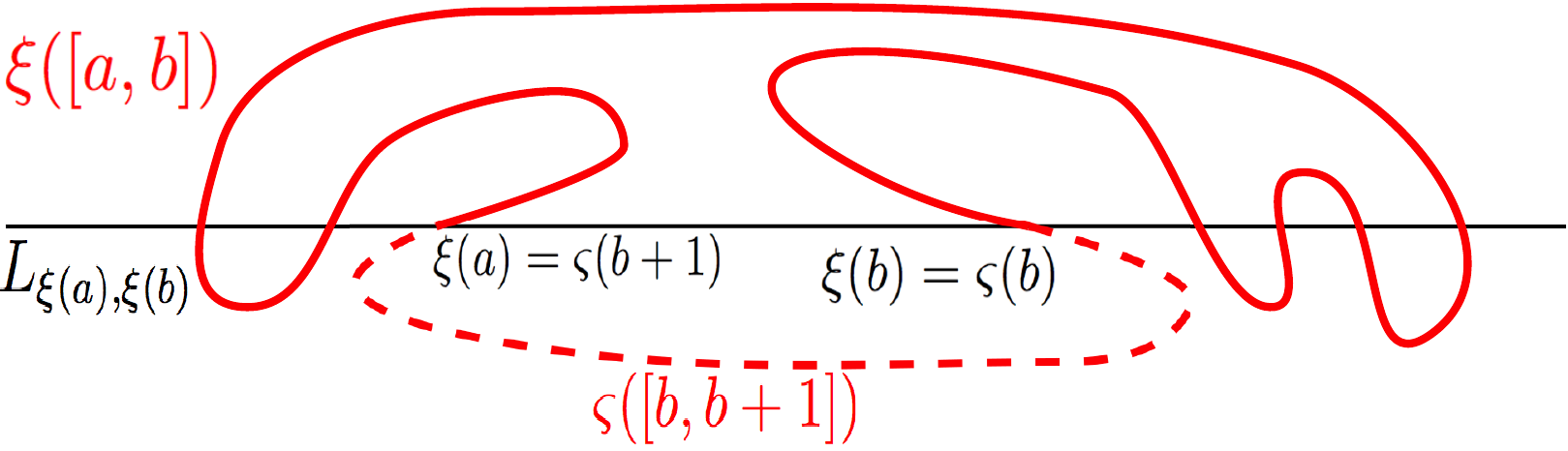}}
\vskip 0.5cm
\subfigure{\includegraphics[scale=0.9]{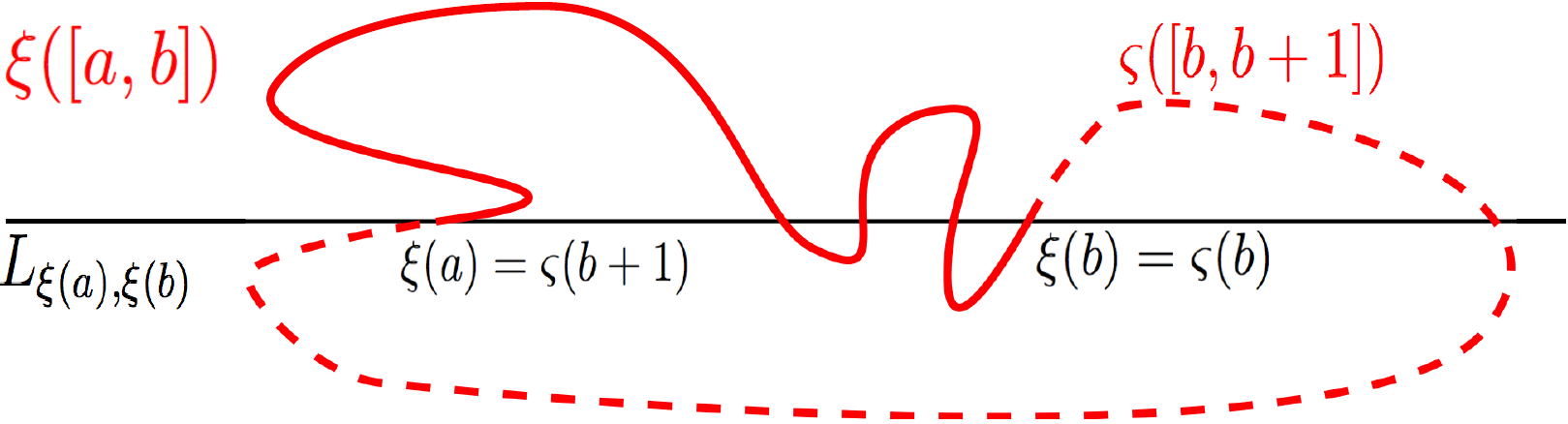}}
\caption{Non-intersecting arcs $\xi([a,b])$ with $\xi([a,b])\cap(\xi(a),\xi(b))=\emptyset$ (top) and $\xi([a,b])\cap\big(L_{\xi(a),\xi(b)}\backslash[\xi(a),\xi(b)]\big)=\emptyset$ (bottom)}
\end{figure}

\begin{lem}\label{nonintersecting}
Let $\xi:\R\to\R^2$ be any $C^1$ embedding, let $\theta$ be one of its continuous arguments, defined as in~\eqref{deftheta}, and let $a<b$ be any two real numbers. If the arc $\xi([a,b])$ is non-intersecting, then $|\theta(a)-\theta(b)|\le4\pi$.
\end{lem}

\noindent{\bf{Proof.}} Let us first consider the case where $\xi([a,b])\cap(\xi(a),\xi(b))=\emptyset$. There is then a $C^1([b,b+1])$ embedding $\varsigma:[b,b+1]\to\R^2$ such that
\beq\label{varsigma1}\left\{\baa{l}
\varsigma(b)=\xi(b),\ \ \varsigma(b+1)=\xi(a),\vspace{3pt}\\
\dot\varsigma(b)=\dot\xi(b),\ \ \dot\varsigma(b+1)=\dot\xi(a),\vspace{3pt}\\
\varsigma([b,b+1])\cap\xi([a,b])=\{\xi(a),\xi(b)\}\eaa\right.
\eeq
and
\beq\label{varsigma2}
|\vartheta(b+1)-\vartheta(b)|\le2\pi,
\eeq
where $\vartheta:[b,b+1]\to\R$ is the unique continuous function such that
\beq\label{varsigma3}
\frac{\dot\varsigma(t)}{|\dot\varsigma(t)|}=\big(\!\cos\vartheta(t),\sin\vartheta(t)\big)\hbox{ for all }t\in[b,b+1]\ \hbox{ and }\ \vartheta(b)=\theta(b)
\eeq
(see Figure~2). Let now $\overline{\xi}$ and $\overline{\theta}$ be the functions defined by
$$\overline{\xi}(t)=\left\{\baa{ll}\xi(t) & \hbox{if }t\in[a,b],\vspace{3pt}\\
\varsigma(t) & \hbox{if }t\in(b,b+1],\eaa\right.\ \hbox{ and }\ \overline{\theta}(t)=\left\{\baa{ll}\theta(t) & \hbox{if }t\in[a,b],\vspace{3pt}\\
\vartheta(t) & \hbox{if }t\in(b,b+1].\eaa\right.$$
These functions are respectively $C^1$ and continuous in $[a,b+1]$. Furthermore,
$$\frac{\dot{\overline{\xi}}(t)}{|\dot{\overline{\xi}}(t)|}=\big(\!\cos\overline{\theta}(t),\sin\overline{\theta}(t)\big)\ \hbox{ for all }t\in[a,b+1]$$
and the closed curve $\overline{\xi}([a,b+1])$ is the boundary of a non-empty domain in $\R^2$ by the third property in~\eqref{varsigma1}. Lastly, $\overline{\xi}$ is one-to-one on $[a,b+1)$, while $\overline{\xi}(b+1)=\overline{\xi}(a)$ and $\dot{\overline{\xi}}(b+1)=\dot{\overline{\xi}}(a)$. Therefore, the continuous argument $\overline{\theta}$ of $\dot{\overline{\xi}}$ is such that
$$|\overline{\theta}(b+1)-\overline{\theta}(a)|=2\pi.$$
In other words, $|\vartheta(b+1)-\theta(a)|=2\pi$. Since $|\vartheta(b+1)-\vartheta(b)|\le2\pi$ and $\vartheta(b)=\theta(b)$ by~\eqref{varsigma2} and~\eqref{varsigma3}, one concludes that
\beq\label{difftheta}
|\theta(b)-\theta(a)|\le4\pi.
\eeq\par
Let us finally consider the case where $\xi([a,b])\cap\big(L_{\xi(a),\xi(b)}\backslash[\xi(a),\xi(b)]\big)=\emptyset$. It is immediate to see that there still exists a $C^1([b,b+1])$ embedding $\varsigma:[b,b+1]\to\R^2$ satisfying~\eqref{varsigma1},~\eqref{varsigma2} and~\eqref{varsigma3}. The remaining arguments are the same as above and~\eqref{difftheta} still holds. The proof of Lemma~\ref{nonintersecting} is thereby complete.\hfill$\Box$\break

Consider again any $C^1$ embedding $\xi:\R\to\R^2$ with continuous argument $\theta$ defined by~\eqref{deftheta} in $\R$. Let $a$ and $b$ be any real numbers such that $a<b$. Denote
$$E=\big\{t\in[a,b];\ \xi(t)\in L_{\xi(a),\xi(b)}\big\}\ \hbox{ and }\ F=[a,b]\backslash E.$$
One has $F=(a,b)\backslash E$ and, by continuity of $\xi$, the set $F$ is an open subset of $(a,b)$. It can then be written as
$$F=\bigcup_{k\in K}I_k,$$
where $K$ is an at most countable set and the sets $I_k$ are pairwise disjoint open intervals
$$I_k=(t_k,t'_k)$$
with $a\le t_k<t'_k\le b$ for every $k\in K$. For every $k\in K$, there holds
$$\xi(t_k)\in L_{\xi(a),\xi(b)}\ \hbox{ and }\ \xi(t'_k)\in L_{\xi(a),\xi(b)}$$
(while $\xi(t)\not\in L_{\xi(a),\xi(b)}$ for any $t\in(t_k,t'_k)$). In particular, the arcs $\xi([t_k,t'_k])$ are all non-intersecting, hence
\beq\label{Ik}
|\theta(t_k)-\theta(t'_k)|\le4\pi\ \hbox{ for all }k\in K
\eeq
by Lemma~\ref{nonintersecting}. For $k\in K$, we say that, relatively to the segment $[\xi(a),\xi(b)]$,
$$\hbox{the arc }\xi([t_k,t'_k])\hbox{ is }\left\{\baa{ll}
\hbox{a middle arc} & \hbox{if }[\xi(t_k),\xi(t'_k)]\subset[\xi(a),\xi(b)],\vspace{3pt}\\
\hbox{a left arc} & \hbox{if }\xi(a)\in(\xi(t_k),\xi(t'_k))\hbox{ and }\xi(b)\not\in(\xi(t_k),\xi(t'_k)),\vspace{3pt}\\
\hbox{a right arc} & \hbox{if }\xi(b)\in(\xi(t_k),\xi(t'_k))\hbox{ and }\xi(a)\not\in(\xi(t_k),\xi(t'_k)),\vspace{3pt}\\
\hbox{a double arc} & \hbox{if }[\xi(a),\xi(b)]\subset(\xi(t_k),\xi(t'_k)),\vspace{3pt}\\
\hbox{an exterior arc} & \hbox{if }[\xi(t_k),\xi(t'_k)]\cap(\xi(a),\xi(b))=\emptyset,\eaa\right.$$
see Figure~3. Notice that these five possibilities are the only ones and that they are pairwise distinct. Let $N_l$, $N_r$ and $N_d$ denote the numbers of left, right and double arcs, respectively, that is, 
$$N_l=\#\big\{k\in K;\ \xi([t_k,t'_k])\hbox{ is a left arc}\big\},$$
and similarly for $N_r$ and $N_d$. These numbers are nonnegative integers or may a priori be $+\infty$.\par
\begin{figure}\centering
\includegraphics[scale=0.9]{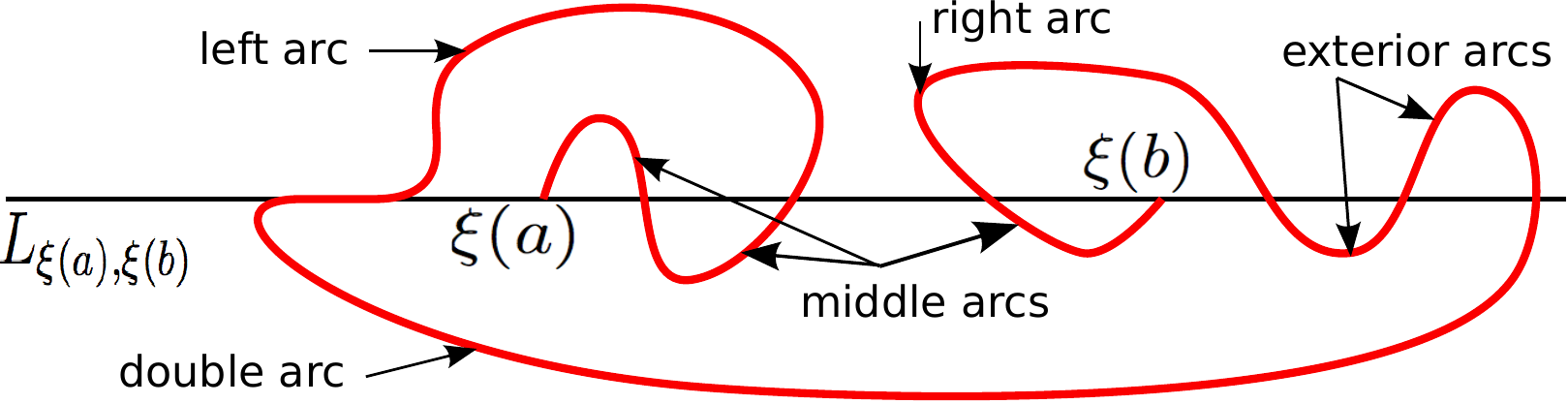}
\caption{Middle, left, right, double and exterior arcs, relatively to the segment $[\xi(a),\xi(b)]$}
\end{figure}
The next lemma shows that these numbers are actually finite and that their sum controls the difference of the argument $\theta$ of $\dot\xi$ between $a$ and $b$.

\begin{lem}\label{arcs}
Under the above notations, the numbers $N_l$, $N_r$ and $N_d$ are nonnegative integers, and
$$|\theta(a)-\theta(b)|\le16\pi\big(N_l+N_r+N_d\big)+4\pi.$$
\end{lem}

\noindent{\bf{Proof.}} Let us first show that the numbers $N_l$, $N_r$ and $N_d$ are finite. Assume by contradiction that, say, $N_l$ is infinite. Then there is a one-to-one map $\rho:\N\to K$ such that $\xi([t_{\rho(n)},t'_{\rho(n)}])$ is a left arc for each $n\in\N$. Since the intervals $(t_{\rho(n)},t'_{\rho(n)})$ are pairwise disjoint and are all included in $[a,b]$, it follows that their length converges to $0$ as $n\to+\infty$, that is, $t'_{\rho(n)}-t_{\rho(n)}\to0$ as $n\to+\infty$. Thus, up to extraction of a subsequence, $t_{\rho(n)}\to T$ and $t'_{\rho(n)}\to T$ as $n\to+\infty$, for some $T\in[a,b]$. Therefore, $\xi(t_{\rho(n)})\to\xi(T)$ and $\xi(t'_{\rho(n)})\to\xi(T)$ as $n\to+\infty$. Since $\xi(a)\in(\xi(t_{\rho(n)}),\xi(t'_{\rho(n)}))$ for each $n\in\N$, one gets that $\xi(a)=\xi(T)$, hence $T=a$ since $\xi$ is one-to-one. Notice also that all real numbers $t_{\rho(n)}$ and $t'_{\rho(n)}$ are different from $a$, since $\xi(a)\in(\xi(t_{\rho(n)}),\xi(t'_{\rho(n)}))$, and that
$$\frac{\xi(t_{\rho(n)})-\xi(a)}{t_{\rho(n)}-a}\to\dot\xi(a)\ \hbox{ and }\ \frac{\xi(t'_{\rho(n)})-\xi(a)}{t'_{\rho(n)}-a}\to\dot\xi(a)\ \hbox{ as }n\to+\infty$$
But the vectors $\xi(t_{\rho(n)})-\xi(a)$ and $\xi(t'_{\rho(n)})-\xi(a)$ point in opposite directions, whereas $a<\min(t_{\rho(n)},t'_{\rho(n)})$ for all $n\in\N$ and $|\dot\xi(a)|>0$. One then gets a contradiction. Therefore, $N_l$ is finite. Similarly, one can show that the number $N_r$ of right arcs is finite.\par
As far as the double arcs are concerned, since the length of each of them is not smaller than $|\xi(b)-\xi(a)|$ and since $\xi$ is one-to-one, one immediately infers that the number $N_d$ of double arcs is finite and
$$N_d\le\frac{{\length}(\xi([a,b]))}{|\xi(b)-\xi(a)|}.$$
As a consequence, the number
$$N=N_l+N_r+N_d$$
is a nonnegative integer.\par
If $N=0$, then the family of arcs $(\xi([t_k,t'_k]))_{k\in K}$ does not contain any left, right or double arcs, and therefore does not contain any exterior arc either. Thus, all arcs $\xi([t_k,t'_k])$ are middle arcs and $\xi(E)\subset[\xi(a),\xi(b)]$. In particular, the arc $\xi([a,b])$ is non-intersecting (the second alternative of~\eqref{defnoninter} is fulfilled) and
\beq\label{4pi0}
|\theta(a)-\theta(b)|\le4\pi\ \hbox{(if }N=0)
\eeq
by Lemma~\ref{nonintersecting}.\par
Let us then assume in the remaining part of the proof that $N\ge 1$. There exist then $N$ non-empty and pairwise distinct intervals $(\tau_i,\tau'_i)\subset(a,b)$ (for $1\le i\le N$) such that, for each $1\le i\le N$, there is a unique $k(i)\in K$ with
\beq\label{tautau'}
(\tau_i,\tau'_i)=I_{k(i)}=(t_{k(i)},t'_{k(i)})
\eeq
and the arc $\xi([\tau_i,\tau'_i])$ is either a left, a right or a double arc. Up to reordering, one can assume without loss of generality that
$$a\le\tau_1<\tau'_1\le\tau_2<\tau'_2\le\cdots\le\tau_N<\tau'_N\le b.$$
Notice that~\eqref{Ik} and~\eqref{tautau'} imply that, for each $1\le i\le N$,
\beq\label{4pi}
|\theta(\tau_i)-\theta(\tau'_i)|\le4\pi.
\eeq\par
For each $1\le i\le N-1$ (if $N\ge2$), the arc $\xi([\tau'_i,\tau_{i+1}])$ does not contain any left, right or double arc and it does not contain $\xi(a)$ or $\xi(b)$ either since $[\tau'_i,\tau_{i+1}]\subset(a,b)$. Hence $\xi(a)\not\in[\xi(\tau'_i),\xi(\tau_{i+1})]$ and $\xi(b)\not\in[\xi(\tau'_i),\xi(\tau_{i+1})]$. As a consequence,
$$\hbox{either }\ [\xi(\tau'_i),\xi(\tau_{i+1})]\subset(\xi(a),\xi(b))\ \hbox{ or }\ [\xi(\tau'_i),\xi(\tau_{i+1})]\cap[\xi(a),\xi(b)]=\emptyset.$$\par
Consider firstly the case where
$$[\xi(\tau'_i),\xi(\tau_{i+1})]\subset(\xi(a),\xi(b))$$
(notice that the relative positions of $\xi(\tau'_i)$ and $\xi(\tau_{i+1})$ between $\xi(a)$ and $\xi(b)$ will not make any difference in the following arguments). Since the continuous curve $\xi([\tau'_i,\tau_{i+1}])$ does not contain any left, right or double arc and it does not contain $\xi(a)$ or $\xi(b)$, it follows that, for any $I_k\subset(\tau'_i,\tau_{i+1})$, the arc $\xi(\overline{I_k})$ is a middle arc (notice that it may well happen that there is no such $I_k$, in which case $\xi([\tau'_i,\tau_{i+1}])=[\xi(\tau'_i),\xi(\tau_{i+1})]$). Furthermore, the arcs $\xi([\tau_i,\tau'_i])$ and $\xi([\tau_{i+1},\tau'_{i+1}])$ cannot be double arcs, and hence they are either left or right arcs. There are four possibilities: left-left, left-right, right-left or right-right. If both are left arcs or both right arcs, then
$$\xi([\tau_i,\tau'_{i+1}])\cap(\xi(\tau_i),\xi(\tau'_{i+1}))=\emptyset$$
and $\xi([\tau_i,\tau'_{i+1}])$ is non-intersecting (see Figure~4, where both arcs $\xi([\tau_i,\tau'_i])$ and $\xi([\tau_{i+1},\tau'_{i+1}])$ are left arcs). If one of the arcs $\xi([\tau_i,\tau'_i])$ and $\xi([\tau_{i+1},\tau'_{i+1}])$ is a left arc and the other one a right arc, then
\begin{figure}\centering
\includegraphics[scale=0.47]{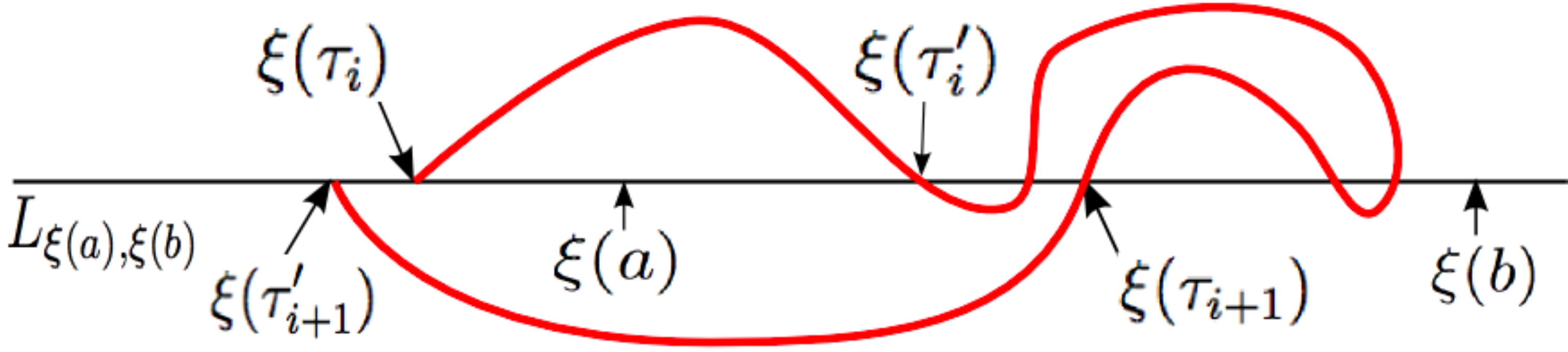}
\caption{The arcs $\xi([\tau_i,\tau'_i])$ and $\xi([\tau_{i+1},\tau'_{i+1}])$ are left arcs, the arc $\xi([\tau_i,\tau'_{i+1}])$ is non-intersecting}
\end{figure}
$$\xi([\tau_i,\tau'_{i+1}])\cap\big(L_{\xi(\tau_i),\xi(\tau'_{i+1})}\backslash[\xi(\tau_i),\xi(\tau'_{i+1})]\big)=\emptyset$$
and $\xi([\tau_i,\tau'_{i+1}])$ is non-intersecting too. Hence, in all cases, $|\theta(\tau_i)-\theta(\tau'_{i+1})|\le4\pi$ by Lemma~\ref{nonintersecting}. But since $|\theta(\tau_i)-\theta(\tau'_i)|\le4\pi$ and $|\theta(\tau_{i+1})-\theta(\tau'_{i+1})|\le4\pi$ by~\eqref{4pi}, one infers that
\beq\label{12pi}
|\theta(\tau'_i)-\theta(\tau_{i+1})|\le12\pi.
\eeq\par
Consider secondly the case where
$$[\xi(\tau'_i),\xi(\tau_{i+1})]\cap[\xi(a),\xi(b)]=\emptyset.$$
Let us then assume without loss of generality that both $\xi(\tau'_i)$ and $\xi(\tau_{i+1})$ are on the left of $\xi(a)$ with respect to $\xi(b)$ (the other case, where $\xi(b)\in(\xi(a),\xi(\tau'_i))\cap(\xi(a),\xi(\tau_{i+1}))$, can be handled similarly). As in the previous paragraph, since the continuous curve $\xi([\tau'_i,\tau_{i+1}])$ does not contain any left or double arc and it does not contain $\xi(a)$, it follows that, for any $I_k=(t_k,t'_k)\subset(\tau'_i,\tau_{i+1})$, the arc $\xi(\overline{I_k})$ is an exterior arc such that both $\xi(t_k)$ and $\xi(t'_k)$ are on the left of $\xi(a)$ with respect to $\xi(b)$. Furthermore, the arcs $\xi([\tau_i,\tau'_i])$ and $\xi([\tau_{i+1},\tau'_{i+1}])$ are either left or double arcs (see Figure~5 in the case where, say, $\xi([\tau_i,\tau'_i])$ is a left arc and $\xi([\tau_{i+1},\tau'_{i+1}])$ is a double arc). It is immediate to check that, in all cases,
\begin{figure}\centering
\includegraphics[scale=0.9]{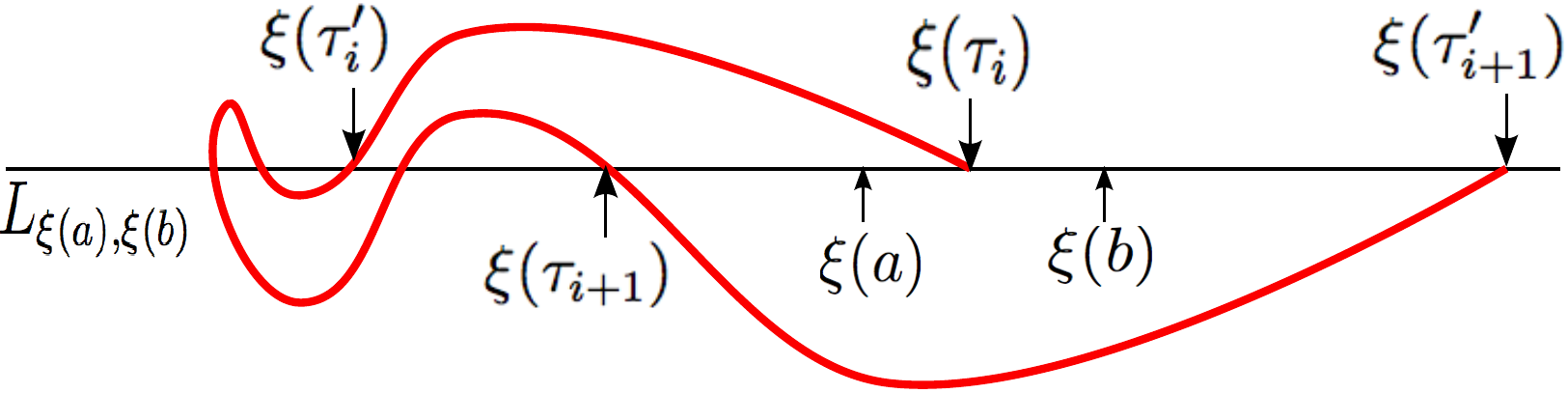}
\caption{The arc $\xi([\tau_i,\tau'_i])$ is a left arc, the arc $\xi([\tau_{i+1},\tau'_{i+1}])$ is a double arc, the arc $\xi([\tau_i,\tau'_{i+1}])$ is non-intersecting}
\end{figure}
$$\xi([\tau_i,\tau'_{i+1}])\cap(\xi(\tau_i),\xi(\tau'_{i+1}))=\emptyset,$$
hence the arc $\xi([\tau_i,\tau'_{i+1}])$ is non-intersecting. Therefore, $|\theta(\tau_i)-\theta(\tau'_{i+1})|\le4\pi$ by Lemma~\ref{nonintersecting} and one infers as in the previous paragraph that
\beq\label{12pibis}
|\theta(\tau'_i)-\theta(\tau_{i+1})|\le12\pi.
\eeq\par
Let us now consider the arc $\xi([a,\tau_1])$, assuming for the moment that $\tau_1>a$. Three cases might occur: either $\xi(\tau_1)$ is on the left of $\xi(a)$ with respect to $\xi(b)$, or $\xi(\tau_1)\in(\xi(a),\xi(b)]$, or $\xi(b)\in(\xi(a),\xi(\tau_1))$. Notice that the latter case is impossible, otherwise $\xi([a,\tau_1])$ would contain a right or double arc, contradicting the definition of $\tau_1$. Thus, either $\xi(\tau_1)$ is on the left of $\xi(a)$ with respect to $\xi(b)$, or $\xi(\tau_1)\in(\xi(a),\xi(b)]$. If $\xi(\tau_1)$ is on the left of $\xi(a)$, since the continuous curve $\xi([a,\tau_1])$ does not contain any left, right or double arc, it follows that, for any $I_k=(t_k,t'_k)\subset(a,\tau_1)$, the arc $\xi(\overline{I_k})$ is an exterior arc such that $\xi(t'_k)$ is on the left of $\xi(a)$ with respect to $\xi(b)$, and either $\xi(t_k)$ is on the left of $\xi(a)$ with respect to $\xi(b)$ or $t_k=a$ with $\xi(t_k)=\xi(a)$. Furthermore, the arc $\xi([\tau_1,\tau'_1])$ is either a left or a double arc, hence $\xi(\tau'_1)$ is on the right of $\xi(a)$ with respect to $\xi(b)$ and
$$\xi([a,\tau'_1])\cap(\xi(a),\xi(\tau'_1))=\emptyset.$$
Therefore, the arc $\xi([a,\tau'_1])$ is non-intersecting and $|\theta(a)-\theta(\tau'_1)|\le4\pi$ by Lemma~\ref{nonintersecting}. One then infers from~\eqref{4pi} with $i=1$ that
$$|\theta(a)-\theta(\tau_1)|\le8\pi.$$
Similarly, if $\xi(\tau_1)\in(\xi(a),\xi(b)]$, since the continuous curve $\xi([a,\tau_1])$ does not contain any left or right arc, it follows that, for any $I_k\subset(a,\tau_1)$, the arc $\xi(\overline{I_k})$ is a middle arc. Furthermore, since $\tau_1<\tau'_1\le b$ (hence, $\xi(\tau_1)$ actually belongs to the interval $(\xi(a),\xi(b))$), the arc $\xi([\tau_1,\tau'_1])$ is either a left or a right arc. If it is a left arc, then $\xi(\tau'_1)$ is on the left of $\xi(a)$ with respect to $\xi(b)$ and
$$\xi([a,\tau'_1])\cap(\xi(a),\xi(\tau'_1))=\emptyset,$$
hence the arc $\xi([a,\tau'_1])$ is non-intersecting. If the arc $\xi([\tau_1,\tau'_1])$ is a right arc, then $\xi(b)\in(\xi(a),\xi(\tau'_1))$ and
$$\xi([a,\tau'_1])\cap\big(L_{\xi(a),\xi(\tau'_1)}\backslash[\xi(a),\xi(\tau'_1)]\big)=\emptyset,$$
hence the arc $\xi([a,\tau'_1])$ is non-intersecting as well. As a consequence, in all cases, Lemma~\ref{nonintersecting} yields $|\theta(a)-\theta(\tau'_1)|\le4\pi$ and $|\theta(\tau_1)-\theta(\tau'_1)|\le4\pi$. Together with~\eqref{4pi} with $i=1$, one gets that
$$|\theta(a)-\theta(\tau_1)|\le8\pi.$$
Notice also that this inequality holds trivially if $\tau_1=a$.\par
Similarly, one can show that
$$|\theta(\tau'_N)-\theta(b)|\le8\pi.$$
As a conclusion, the previous two inequalities together with~\eqref{4pi},~\eqref{12pi} and~\eqref{12pibis} imply that
$$\baa{l}
|\theta(a)-\theta(b)|\vspace{3pt}\\
\qquad\le\displaystyle|\theta(a)-\theta(\tau_1)|+\sum_{i=1}^{N-1}\big(|\theta(\tau_i)-\theta(\tau'_i)|+|\theta(\tau'_i)-\theta(\tau_{i+1})|\big)+|\theta(\tau_N)-\theta(\tau'_N)|+|\theta(\tau'_N)-\theta(b)|\vspace{3pt}\\
\qquad\le8\pi+(N-1)(4\pi+12\pi)+4\pi+8\pi=16\pi N+4\pi,\eaa$$
which is the desired result. Reminding~\eqref{4pi0}, the proof of Lemma~\ref{arcs} is thereby complete.\hfill$\Box$


\subsection{Logarithmic growth of the argument of $\nabla u$ along its trajectories}\label{sec32}

For any $C^2(\R^2)$ solution $w$ of~\eqref{1} satisfying~\eqref{etaw} with $0<\eta\le1$, let us remind the definitions~\eqref{argv} and~\eqref{defvarphi} of the arguments $\phi_w$ and $\varphi_w$ of $w$ and $\nabla u_w$, and that of the parametrizations $\sigma:\R\to\R^2$, defined by $\dot\sigma(t)=\nabla u_w(\sigma(t))$, of the trajectories $\Sigma$ of the gradient flow. We also recall that ${\osc}_E\,\phi_w={\osc}_E\,\varphi_w$ for any non-empty subset $E\subset\R^2$.\par
The present section is devoted to the proof of some estimates on the logarithmic growth of the argument $\varphi_w$ of $\nabla u_w$ along its trajectories. For the sake of simplicity, we drop in the sequel the indices $w$, that is, we write $u=u_w$, $\phi=\phi_w$ and $\varphi=\varphi_w$. In order to state these estimates, we first show an auxiliary elementary lemma on the local behavior of the trajectories of the gradient flow around any point.

\begin{lem}\label{oneleft}
For any $\eta\in(0,1]$, for any $C^2(\R^2)$ solution $w$ of~\eqref{1} satisfying~\eqref{etaw} and for any trajectory of the gradient flow, with parametrization $\dot\sigma(t)=\nabla u(\sigma(t))$ for $t\in\R$, there are no real numbers $(\tau_i)_{1\le i\le 4}$ satisfying
\beq\label{tau14}
\tau_1<\tau_2<\tau_3\le\tau_4\ \hbox{ or }\ \tau_1>\tau_2>\tau_3\ge\tau_4,
\eeq
and such that
\beq\label{xitau14}
\mathop{\osc}_{B(\sigma(\tau_1),1)}\phi<\frac{\pi}{2},\ \ \sigma(\tau_1)\in(\sigma(\tau_2),\sigma(\tau_3)),\ \hbox{ and }\ |\sigma(\tau_1)-\sigma(\tau_4)|<\eta^4.
\eeq 
\end{lem}

\noindent{\bf{Proof.}} Assume by way of contradiction that there exist a trajectory of the gradient flow, with parametrization $\dot\sigma(t)=\nabla u(\sigma(t))$ for $t\in\R$, and some real numbers $(\tau_i)_{1\le i\le 4}$ satisfying~\eqref{tau14} and~\eqref{xitau14}. Let us only consider the first case
$$\tau_1<\tau_2<\tau_3\le\tau_4$$
in~\eqref{tau14} (the second case $\tau_1>\tau_2>\tau_3\ge\tau_4$ can be obtained from the first case by replacing~$w$ by~$-w$ and~$u$ by~$-u$). The mean value theorem together with~\eqref{etaw} and~\eqref{xitau14} then yields 
$$|u(\sigma(\tau_1))-u(\sigma(\tau_4))|<\eta^3.$$
Actually, since the function $t\mapsto g(t):=u(\sigma(t))$ is increasing in $\R$ ($g'(t)=|\nabla u(\sigma(t))|^2\ge\eta^2>0$ in~$\R$), one gets that $u(\sigma(\tau_1))<u(\sigma(\tau_2))<u(\sigma(\tau_3))\le u(\sigma(\tau_4))$, and finally
$$u(\sigma(\tau_1))<u(\sigma(\tau_2))<u(\sigma(\tau_3))\le u(\sigma(\tau_4))<u(\sigma(\tau_1))+\eta^3.$$
Since $g'(t)\ge\eta^2>0$ for all $t\in\R$, this yields
$$\tau_1<\tau_2<\tau_3\le\tau_4<\tau_1+\eta.$$\par
Now, for every $\tau\in[\tau_1,\tau_4]$, there holds $|\dot\sigma(\tau)|=|\nabla u(\sigma(\tau))|\le\eta^{-1}$, hence $|\sigma(\tau)-\sigma(\tau_1)|\le\eta^{-1}(\tau_4-\tau_1)<1$. Since ${\osc}_{B(\sigma(\tau_1),1)}\varphi={\osc}_{B(\sigma(\tau_1),1)}\phi<\pi/2$ by assumption~\eqref{xitau14} and since $\varphi$ is the argument of $\nabla u$, one then gets that $\nabla u(\sigma(\tau))\cdot\nabla u(\sigma(\tau_1))>0$ for all $\tau\in[\tau_1,\tau_4]$, that is, $\dot\sigma(\tau)\cdot\dot\sigma(\tau_1)>0$. This implies that $(\sigma(\tau)-\sigma(\tau_1))\cdot\dot\sigma(\tau_1)>0$ for all $\tau\in(\tau_1,\tau_4]$. In particular,
$$(\sigma(\tau_2)-\sigma(\tau_1))\cdot\dot\sigma(\tau_1)>0\ \hbox{ and }\ (\sigma(\tau_3)-\sigma(\tau_1))\cdot\dot\sigma(\tau_1)>0.$$
But the vectors $\sigma(\tau_2)-\sigma(\tau_1)$ and $\sigma(\tau_3)-\sigma(\tau_1)$ point in opposite directions, since $\sigma(\tau_1)\in(\sigma(\tau_2),\sigma(\tau_3))$ by assumption~\eqref{xitau14}. This leads to a contradiction and the proof of Lemma~\ref{oneleft} is thereby complete.\hfill$\Box$\break

In order to state the key Lemma~\ref{log} below on the logarithmic growth of the argument $\varphi$ of $\nabla u$ along its trajectories, let us introduce a few auxiliary notations. For $\eta\in(0,1]$, denote
\beq\label{defC1eta}
C_1(\eta):=\sup_{t\in[1,+\infty)}h_1(t),\ \hbox{ where }h_1(t)=288\pi\eta^{-2}\times\frac{\ln_2(t\eta^{-4})+2}{\ln(3+t)}
\eeq
and $\ln_2t=\ln t/\ln 2$ for any $t>0$. Notice immediately that $C_1(\eta)\in\R$ and
\beq\label{C2bis}
C_1(\eta)\ge\lim_{t\to+\infty}h_1(t)=\frac{288\pi\eta^{-2}}{\ln 2}>96\pi\eta^{-2}\ge96\pi>\frac{\pi}{2}.
\eeq
With this constant $C_1(\eta)>0$, the following estimate holds for the trajectories of the gradient flow.

\begin{lem}\label{log}
For any $\eta\in(0,1]$, for any $C^2(\R^2)$ solution $w$ of~\eqref{1} satisfying~\eqref{etaw} and for any $x\in\R^2$, there holds
$$\mathop{\osc}_{B(x,1)}\varphi=\mathop{\osc}_{B(x,1)}\phi<\frac{\pi}{2}\ \Longrightarrow\ \Big(\forall\,y\in\Sigma_x,\ \ |\varphi(x)-\varphi(y)|\le C_1(\eta)\,\ln\!\big(3+|x-y|\big)\Big).$$
\end{lem}

The general underlying idea of the proof is the following: if a trajectory of the gradient flow turns many times around some points then the stream function $u$ would become large along it and the oscillations of $u$ between some not-too-far points would be large. This would lead to a contradiction, since $u$ is Lipschitz continuous. The detailed proof of Lemma~\ref{log} is actually much more involved, since one needs to treat the cases of left, right or double arcs around some segments.\hfill\break

\noindent{\bf{Proof.}} Assume by contradiction that the conclusion of Lemma~\ref{log} does not hold for some $C^2(\R^2)$ solution $w$ of~\eqref{1} satisfying~\eqref{etaw} with $0<\eta\le1$. Then, there exist $x\in\R^2$ and $y\in\Sigma_x$ such that ${\osc}_{B(x,1)}\varphi<\pi/2$ and
$$|\varphi(x)-\varphi(y)|>C_1(\eta)\,\ln\!\big(3+|x-y|\big).$$
Let $\sigma$ be a parametrization of $\Sigma_x$ solving $\dot\sigma(t)=\nabla u(\sigma(t))$ for $t\in\R$, and $a,b\in\R$ be such that
$$x=\sigma(a)\ \hbox{ and }\ y=\sigma(b).$$
One has $|\varphi(x)-\varphi(y)|>C_1(\eta)>\pi/2$ (by~\eqref{C2bis}), and the property ${\osc}_{B(x,1)}\varphi<\pi/2$ implies that
\beq\label{xiab}
|\sigma(a)-\sigma(b)|=|x-y|\ge1.
\eeq
Let us assume here that $a<b$ (the case $a>b$ can be handled similarly). Let $\theta$ be the continuous argument of $\dot\sigma$, defined as in~\eqref{deftheta} with the embedding $\xi:=\sigma$. Since $\dot\sigma(t)=\nabla u(\sigma(t))$ in $\R$, it follows from~\eqref{defvarphi} and the continuity of $\theta$ and $\varphi(\sigma(\cdot))$ in $\R$ that there is an integer $q\in\Z$ such that
$$\theta(t)=\varphi(\sigma(t))+2\pi q\ \hbox{ for all }t\in\R.$$
In particular, $\theta(a)-\theta(b)=\varphi(\sigma(a))-\varphi(\sigma(b))=\varphi(x)-\varphi(y)$, hence
\beq\label{thetaxi}
|\theta(a)-\theta(b)|>C_1(\eta)\,\ln\big(3+|x-y|\big)=C_1(\eta)\,\ln\big(3+|\sigma(a)-\sigma(b)|\big)
>C_1(\eta)>96\pi
\eeq
by~\eqref{C2bis}.\par
With the same notations as in Section~\ref{sec31}, let $N_l$, $N_r$ and $N_d$ be the (finite) numbers of left, right and double arcs contained in the curve $\sigma([a,b])$, relatively to the segment $[\sigma(a),\sigma(b)]$. It follows from the inequalities $|\theta(a)-\theta(b)|>96\pi>8\pi$ and Lemma~\ref{arcs} that
\beq\label{theta16}
|\theta(a)-\theta(b)|\le16\pi(N_l+N_r+N_d)+4\pi<16\pi(N_l+N_r+N_d)+\frac{|\theta(a)-\theta(b)|}{2},
\eeq
hence $|\theta(a)-\theta(b)|<32\pi(N_l+N_r+N_d)$ and
\beq\label{theta96}
\max\big(N_l,N_r,N_d\big)>\frac{|\theta(a)-\theta(b)|}{96\pi}.
\eeq
Therefore, at least one of the three numbers $N_l$, $N_r$ or $N_d$ is larger than $|\theta(a)-\theta(b)|/(96\pi)$. The three cases will be treated separately.\par
\vskip 0.2cm
{\it Case 1: $N_l>|\theta(a)-\theta(b)|/(96\pi)$.} Notice that~\eqref{thetaxi} yields $N_l>1$, that is, $N_l\ge2$. By definition of $N_l$ and of a left arc, there are some real numbers
$$a<\rho_1<\rho'_1\le\rho_2<\rho'_2\le\cdots\le\rho_{N_l}<\rho'_{N_l}\le b$$
such that, for each $1\le i\le N_l$, there is $k(i)\in K$ such that $(\rho_i,\rho'_i)=I_{k(i)}=(t_{k(i)},t'_{k(i)})$ and $\sigma([\rho_i,\rho'_i])$ is a left arc (we still use the same notations $E$, $F$ and $(I_k)_{k\in K}$ for the arc $\sigma([a,b])$ as in Section~\ref{sec31} for $\xi([a,b])$). We also recall that, for any left arc $\sigma([\rho_i,\rho'_i])$, one has $\sigma(a)\in(\sigma(\rho_i),\sigma(\rho'_i))$, hence $\rho_i$ and $\rho'_i$ are (strictly) larger than $a$.\par
Define
\beq\label{defm}
m=\big[\ln_2(|\sigma(a)-\sigma(b)|\,\eta^{-4})\big]+2,
\eeq
where, for any $x\in\R$, $[x]$ denotes the integer part of $x$. Notice that~\eqref{xiab} together with $0<\eta\le1$ implies that $m$ is an integer such that
$$m\ge2.$$
Divide now the segment $[\sigma(a),\sigma(b)]$ into $m$ segments $[x_p,y_p]$ (for $p=1,\cdots,m$) defined by
$$\left\{\baa{ll}
\displaystyle x_1=\frac{\sigma(a)+\sigma(b)}{2},\ y_1=\sigma(b), & \ \ \displaystyle x_2=\frac{\sigma(a)+x_1}{2},\ y_2=x_1,\vspace{3pt}\\
\cdots & \\
\displaystyle x_{m-1}=\frac{\sigma(a)+x_{m-2}}{2},\ y_{m-1}=x_{m-2}, & \ \ \displaystyle x_m=\sigma(a)\hbox{ and }y_m=x_{m-1},\eaa\right.$$
see Figure~6. It is immediate to see that the sets $(x_p,y_p]$ form a partition of $(\sigma(a),\sigma(b)]$, that
\begin{figure}\centering
\includegraphics[scale=0.8]{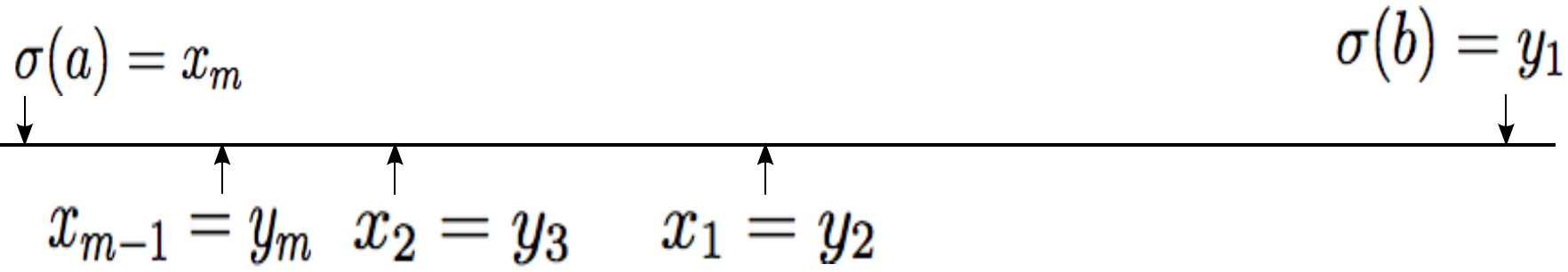}
\caption{The points $\sigma(a)=x_m$, $x_{m-1}=y_m,\ldots,x_1=y_2$, $\sigma(b)=y_1$ in case $m=4$}
\end{figure}
$$|x_p-y_p|=2^{-p}|\sigma(a)-\sigma(b)|\hbox{ for each }1\le p\le m-1,\ \ |x_m-y_m|=2^{1-m}|\sigma(a)-\sigma(b)|$$
and that
\beq\label{xpyp}
|\sigma(a)-x_p|=2^{-p}|\sigma(a)-\sigma(b)|\hbox{ and }|\sigma(a)-y_p|=2^{1-p}|\sigma(a)-\sigma(b)|,\hbox{ for each }1\le p\le  m-1.
\eeq
By definition of $m$ in~\eqref{defm}, it also follows that
\beq\label{eta4}
|\sigma(a)-y_m|=|x_m-y_m|=2^{1-m}|\sigma(a)-\sigma(b)|<\eta^4.
\eeq\par
Coming back to the family of the left arcs $(\sigma([\rho_i,\rho'_i]))_{1\le i\le N_l}$, we remind that, by definition of a left arc, either $\sigma(\rho_i)$ or $\sigma(\rho'_i)$ belongs to $(\sigma(a),\sigma(b)]$, for any $1\le i\le N_l$. Since the sets $(x_p,y_p]$ with $1\le p\le m$ are a partition of $(\sigma(a),\sigma(b)]$, it follows that, if we denote, for $1\le p\le m$,
$$J_p:=\big\{i\in\{1,\ldots,N_l\};\ \sigma(\rho_i)\in(x_p,y_p]\hbox{ or }\sigma(\rho'_i)\in(x_p,y_p]\big\},$$
then
\beq\label{Jp}
\bigcup_{1\le p\le m}J_p=\{1,\ldots,N_l\}.
\eeq\par
We now claim that
\beq\label{leftarc}
\#J_m\le 1,
\eeq
that is, there is at most one left arc intersecting $(x_m,y_m]=(\sigma(a),y_m]$. Assume by contradiction that there are two integers $i$ and $j$ with $1\le i<j\le N_l$ and such that either $\sigma(\rho_i)$ or $\sigma(\rho'_i)$ belongs to $(\sigma(a),y_m]$ and either $\sigma(\rho_j)$ or $\sigma(\rho'_j)$ belongs to $(\sigma(a),y_m]$. Since $|\sigma(a)-y_m|<\eta^4$ by~\eqref{eta4}, one infers in particular that
$$\min\big(|\sigma(\rho_j)-\sigma(a)|,|\sigma(\rho'_j)-\sigma(a)|\big)<\eta^4.$$
Furthermore, $\sigma(a)\in(\sigma(\rho_i),\sigma(\rho'_i))$ since $\sigma([\rho_i,\rho'_i])$ is a left arc. Since
$$\mathop{\osc}_{B(\sigma(a),1)}\phi=\mathop{\osc}_{B(x,1)}\phi=\mathop{\osc}_{B(x,1)}\varphi<\frac{\pi}{2},$$
Lemma~\ref{oneleft} applied with $(\tau_1,\tau_2,\tau_3)=(a,\rho_i,\rho'_i)$ and $\tau_4=\rho_j\,(\ge\rho'_i)$ if $|\sigma(\rho_j)-\sigma(a)|<\eta^4$ (resp. $\tau_4=\rho'_j\,(>\rho'_i)$ if $|\sigma(\rho'_j)-\sigma(a)|<\eta^4$) leads to a contradiction. Therefore, the claim~\eqref{leftarc} has been proved.\par
Putting together~\eqref{Jp} and~\eqref{leftarc}, it follows that there is $p\in\{1,\ldots,m-1\}$ such that
\beq\label{defn}
n:=\#J_p\ge\frac{N_l-1}{m-1}
\eeq
(we also remind that $N_l\ge 2$ and $m\ge 2$, whence $n$ is a positive integer). Write $J_p=\{i_1,\ldots,i_n\}$ with $1\le i_1<\cdots<i_n\le N_l$. Thus,
\beq\label{rhorho'}
a<\rho_{i_1}<\rho'_{i_1}\le\rho_{i_2}<\rho'_{i_2}\le\cdots\le\rho_{i_n}<\rho'_{i_n}\le b.
\eeq
For each $1\le j\le n$, one has $\sigma(a)\in(\sigma(\rho_{i_j}),\sigma(\rho'_{i_j}))$ (because $\sigma([\rho_{i_j},\rho'_{i_j}])$ is a left arc), while either $\sigma(\rho_{i_j})$ or $\sigma(\rho_{i_j})$ belongs to $(x_p,y_p]\,(\subset(\sigma(a),\sigma(b)])$. Hence, $(\sigma(a),x_p]\subset(\sigma(\rho_{i_j}),\sigma(\rho'_{i_j}))$ and
$${\length}\big(\sigma([\rho_{i_j},\rho'_{i_j}])\big)\ge|\sigma(\rho_{i_j})-\sigma(\rho'_{i_j})|\ge|\sigma(a)-x_p|=2^{-p}|\sigma(a)-\sigma(b)|$$
by~\eqref{xpyp}. Since $\sigma$ is an embedding, one also infers from~\eqref{rhorho'} that
$${\length}\big(\sigma([a,\rho_{i_n}])\big)\ge\sum_{1\le j\le n-1}{\length}\big(\sigma([\rho_{i_j},\rho'_{i_j}])\big)\ge(n-1)\,2^{-p}\,|\sigma(a)-\sigma(b)|.$$
Furthermore, Lemma~\ref{lemsigma} yields
\beq\label{arhoin}
u(\sigma(\rho'_{i_n}))>u(\sigma(\rho_{i_n}))\ge u(\sigma(a))+\eta\times{\length}\big(\sigma([a,\rho_{i_n}])\big),
\eeq
hence
\beq\label{ineq1}
u(\sigma(\rho'_{i_n}))>u(\sigma(\rho_{i_n}))\ge u(\sigma(a))+\eta\,(n-1)\,2^{-p}|\sigma(a)-\sigma(b)|.
\eeq\par
On the other hand, one knows that either $\sigma(\rho_{i_n})$ or $\sigma(\rho'_{i_n})$ belongs to $(x_p,y_p]\,(\subset(\sigma(a),\sigma(b)])$, since $i_n\in J_p$. Thus,
$$\min\big(|\sigma(\rho_{i_n})-\sigma(a)|,|\sigma(\rho'_{i_n})-\sigma(a)|\big)\le|\sigma(a)-y_p|=2^{1-p}|\sigma(a)-\sigma(b)|$$
by~\eqref{xpyp}. Since $|\nabla u(x)|\le\eta^{-1}$ for all $x\in\R^2$, it follows that
$$\min\big(|u(\sigma(\rho_{i_n}))-u(\sigma(a))|,|u(\sigma(\rho'_{i_n}))-u(\sigma(a))|\big)\le 2^{1-p}\eta^{-1}|\sigma(a)-\sigma(b)|.$$
But, as already underlined, $u(\sigma(a))<u(\sigma(\rho_{i_n}))<u(\sigma(\rho'_{i_n}))$. Consequently,
\beq\label{ineq2}
u(\sigma(\rho_{i_n}))-u(\sigma(a))\le2^{1-p}\eta^{-1}|\sigma(a)-\sigma(b)|.
\eeq\par
From~\eqref{ineq1} and~\eqref{ineq2}, one infers that
$$\eta\,(n-1)\,2^{-p}|\sigma(a)-\sigma(b)|\le2^{1-p}\eta^{-1}|\sigma(a)-\sigma(b)|$$
and, since $|\sigma(a)-\sigma(b)|\ge1>0$ by~\eqref{xiab} and $0<\eta\le1$, it follows that
$$n\le2\eta^{-2}+1\le3\eta^{-2}.$$
Together with~\eqref{defn}, with our assumption $N_l>|\theta(a)-\theta(b)|/(96\pi)$ and with $\eta^{-2}\ge1$, one gets in particular that
$$|\theta(a)-\theta(b)|<96\pi\times3\eta^{-2}\times(m-1)+96\pi<288\pi\eta^{-2}\times m.$$
From~\eqref{thetaxi} and the definition~\eqref{defm} of $m$, one infers that
$$C_1(\eta)\,\ln(3+|\sigma(a)-\sigma(b)|)<|\theta(a)-\theta(b)|<288\pi\eta^{-2}\times\big(\big[\ln_2(|\sigma(a)-\sigma(b)|\,\eta^{-4})\big]+2\big),$$
hence
$$C_1(\eta)<288\pi\eta^{-2}\times\frac{\ln_2(|\sigma(a)-\sigma(b)|\,\eta^{-4})+2}{\ln(3+|\sigma(a)-\sigma(b)|)}.$$
Since $|\sigma(a)-\sigma(b)|\ge1$ by~\eqref{xiab}, the previous inequality contradicts the definition~\eqref{defC1eta} of the constant $C_1(\eta)$.\par
\vskip 0.2cm
{\it Case 2: $N_r>|\theta(a)-\theta(b)|/(96\pi)$.} The study of this case is similar to that of Case~1 (one especially uses Lemma~\ref{oneleft} with the second alternative in~\eqref{tau14} to show that the number of right arcs meeting a small portion of $[\sigma(a),\sigma(b)]$ around the point $\sigma(b)$ is at most~$1$), and one similarly gets a contradiction with the definition of $C_1(\eta)$.\par
\vskip 0.2cm
{\it Case 3: $N_d>|\theta(a)-\theta(b)|/(96\pi)$.} For each double arc $\sigma([t_k,t'_k])$, one has by definition $[\sigma(a),\sigma(b)]\subset(\sigma(t_k),\sigma(t'_k))$, hence
$${\length}(\sigma([t_k,t'_k]))\ge|\sigma(a)-\sigma(b)|.$$
Since $\sigma:[a,b]\to\R^2$ is an embedding and $\sigma([a,b])$ contains $N_d$ double arcs, it follows that ${\length}(\sigma([a,b]))\ge N_d\,|\sigma(a)-\sigma(b)|$. Lemma~\ref{lemsigma} then yields
$$u(\sigma(b))-u(\sigma(a))\ge\eta\times{\length}(\sigma([a,b]))\ge \eta\,N_d\,|\sigma(a)-\sigma(b)|.$$
On the other hand, $u(\sigma(b))-u(\sigma(a))\le\eta^{-1}|\sigma(b)-\sigma(a)|$ by~\eqref{etaw} and, since $\sigma(a)\neq\sigma(b)$, one infers that $\eta\,N_d\le\eta^{-1}$, hence
$$|\theta(a)-\theta(b)|<96\pi N_d\le96\pi\eta^{-2}$$
from our assumption in Case 3. Since $|\theta(a)-\theta(b)|>C_1(\eta)\ln\big(3+|\sigma(a)-\sigma(b)|\big)>C_1(\eta)$ by~\eqref{thetaxi}, one finally gets that $C_1(\eta)<96\pi\eta^{-2}$, contradicting~\eqref{C2bis}. Case 3 is thus ruled out too and the proof of Lemma~\ref{log} is thereby complete.\hfill$\Box$


\subsection{Logarithmic growth of the argument of $\nabla u$ along the streamlines}\label{sec33}

For any $C^2(\R^2)$ solution $w$ of~\eqref{1} satisfying~\eqref{etaw} with $0<\eta\le1$, let us remind the definition of the parametrizations $\gamma:\R\to\R^2$, given by $\dot\gamma(t)=w(\gamma(t))$, of the streamlines $\Gamma$ of the flow~$w$. The present section is devoted to the proof of some estimates on the logarithmic growth of the argument $\varphi=\varphi_w$ of $\nabla u=\nabla u_w$ along the streamlines. In order to state these estimates, we first show an auxiliary elementary lemma, a bit similar to Lemma~\ref{oneleft}, on the local behavior of a streamline around a point.

\begin{lem}\label{oneleftbis}
For any $\eta\in(0,1]$, for any $C^2(\R^2)$ solution $w$ of~\eqref{1} satisfying~\eqref{etaw} and for any streamline of the flow, with parametrization $\dot\gamma(t)=w(\gamma(t))$ for $t\in\R$, there are no real numbers $(\tau_i)_{1\le i\le 4}$ satisfying
\beq\label{tau14bis}
\tau_1<\tau_2<\tau_3\le\tau_4\ \hbox{ or }\ \tau_1>\tau_2>\tau_3\ge\tau_4,
\eeq
and such that
\beq\label{xitau14bis}
\mathop{\osc}_{B(\gamma(\tau_1),1)}\varphi=\mathop{\osc}_{B(\gamma(\tau_1),1)}\phi<\frac{\pi}{4},\ \ \gamma(\tau_1)\in(\gamma(\tau_2),\gamma(\tau_3)),\ \hbox{ and }\ |\gamma(\tau_1)-\gamma(\tau_4)|<\frac{\eta^2}{4}.
\eeq 
\end{lem}

\noindent{\bf{Proof.}} The beginning of the proof is similar to that of Lemma~\ref{oneleft}, but the end differs, due to the fact that the function $u$ is not anymore strictly monotone, but constant, along the streamlines. Assume by way of contradiction that there exist a trajectory of the flow, with parametrization $\dot\gamma(t)=w(\gamma(t))$ for $t\in\R$, and some real numbers $(\tau_i)_{1\le i\le 4}$ satisfying~\eqref{tau14bis} and~\eqref{xitau14bis}. Let us consider only the first case
$$\tau_1<\tau_2<\tau_3\le\tau_4$$
in~\eqref{tau14bis} (the second case $\tau_1>\tau_2>\tau_3\ge\tau_4$ can be obtained from the first case by replacing $w$ by $-w$ and $u$ by $-u$). Up to shifting $t$ and translating and rotating the frame, one can assume without loss of generality that
$$\tau_1=0,\ \ \gamma(\tau_1)=\gamma(0)=(0,0)=0,\ \ \hbox{ and }\ w(0)=(|w(0)|,0).$$
In particular, $\gamma=\gamma_0$ and $\big\{\gamma(t);\ t\in\R\big\}=\Gamma_0$.\par
Let $Q$ be the open rectangle centered at $0$ and defined by
$$Q=\Big(-\frac{\eta^2}{4},\frac{\eta^2}{4}\Big)\times\Big(-\frac12,\frac12\Big).$$
Let us show in this paragraph that $\Gamma_0\cap Q$ is a graph in the variable $x_1\in(-\eta^2/4,\eta^2/4)$. First of all, by Lemma~\ref{unbounded}, there are some first exit times $\tau^-$ and $\tau^+$ (from the rectangle $Q$) such that
$$\tau^-<0<\tau^+,\ \ \gamma(\tau^{\pm})\in\partial Q,\ \hbox{ and }\ \gamma(\tau)\in Q\hbox{ for all }\tau\in(\tau^-,\tau^+).$$
Denote $\gamma=(\gamma_1,\gamma_2)$ the two coordinates of $\gamma$. For all $\tau\in(\tau^-,\tau^+)$, there holds $\gamma(\tau)\in Q\subset B(0,1)$ and, since ${\osc}_{B(0,1)}\phi={\osc}_{B(\gamma(\tau_1),1)}\phi<\pi/4$ by assumption~\eqref{xitau14bis}, one has $w_1(\gamma(\tau))\,|w(0)|=w(\gamma(\tau))\cdot w(0)>0$ and $|w_2(\gamma(\tau))|<w_1(\gamma(\tau))$. In other words,
\beq\label{dotgamma}
\dot\gamma_1(\tau)>0\hbox{ and }|\dot\gamma_2(\tau)|<\dot\gamma_1(\tau)\ \hbox{ for all }\tau\in(\tau^-,\tau^+).
\eeq
Since $\eta\le|\dot\gamma(\tau)|=|w(\gamma(\tau))|\le\eta^{-1}$, it follows that $\dot\gamma_1(\tau)>\eta/\sqrt{2}$ and $|\dot\gamma_2(\tau)|<\eta^{-1}/\sqrt{2}$ for all $\tau\in(\tau^-,\tau^+)$. Remembering that $\gamma(0)=(0,0)$ and $|\gamma_1(\tau^{\pm})|\le\eta^2/4$ (since $\gamma(\tau^{\pm})\in\partial Q$), one then gets that
$$-\frac{\eta}{2\sqrt{2}}\le\tau^-<0<\tau^+\le\frac{\eta}{2\sqrt{2}}\ \hbox{ and }\ |\gamma_2(\tau)|<\frac{\eta^{-1}}{\sqrt{2}}\times\frac{\eta}{2\sqrt{2}}=\frac14<\frac12$$
for all $\tau\in(\tau^-,\tau^+)$. Therefore, the exit points $\gamma(\tau^{\pm})$ belong to the lateral sides of the rectangle $Q$, namely $|\gamma_1(\tau^{\pm})|=\eta^2/4$. Moreover, remembering that ${\dot\gamma}_1(\tau)\ge\eta/\sqrt{2}>0$ for all $\tau\in(\tau^-,\tau^+)$, it follows that $\gamma_1(\tau^{\pm})=\pm\eta^2/4$ and that, for each $x_1\in(-\eta^2/4,\eta^2/4)$, there is a unique real number $\tau_{x_1}\in(\tau^-,\tau^+)$ such that $\gamma_1(\tau_{x_1})=x_1$. The curve $\gamma((\tau^-,\tau^+))=\big\{\gamma(\tau);\ \tau\in(\tau^-,\tau^+)\big\}$ can then be written as a graph
$$\gamma((\tau^-,\tau^+))=\Big\{(x_1,\gamma_2(\tau_{x_1}));\ x_1\in\Big(-\frac{\eta^2}{4},\frac{\eta^2}{4}\Big)\Big\}.$$
On the other hand, there holds $\nabla u(0)=(0,-|\nabla u(0)|)$ (as $\nabla^\perp u(0)=w(0)=(|w(0)|,0)$) and ${\osc}_{Q}\varphi\le{\osc}_{B(0,1)}\varphi={\osc}_{B(\gamma(\tau_1),1)}\varphi={\osc}_{B(\gamma(\tau_1),1)}\phi<\pi/4<\pi/2$ by assumption~\eqref{xitau14bis}. Hence,
$$\frac{\partial u}{\partial x_2}<0\ \hbox{in }Q.$$
Finally, since $u$ is constant along $\Gamma_0$ (that is, $t\mapsto u(\gamma(t))$ is constant in $\R$), one concludes that
$$\Gamma_0\cap Q=\gamma((\tau^-,\tau^+))=\Big\{(x_1,\gamma_2(\tau_{x_1}));\ x_1\in\Big(-\frac{\eta^2}{4},\frac{\eta^2}{4}\Big)\Big\}.$$\par
Remember now the assumption~\eqref{xitau14bis} (with $\tau_1=0$ and $\gamma(0)=0$). The inequality $|\gamma(\tau_4)|=|\gamma(\tau_1)-\gamma(\tau_4)|<\eta^2/4\,(<1/2)$ yields $\gamma(\tau_4)\in\Gamma_0\cap Q$, hence $\tau_4\in(\tau^-,\tau^+)$ and finally
$$\tau^-<0=\tau_1<\tau_2<\tau_3\le\tau_4<\tau^+.$$
As a consequence, one has $\dot\gamma_1(\tau)>0$ for all $\tau\in[0,\tau_4]$ by~\eqref{dotgamma}. In particular, $0=\gamma_1(0)=\gamma_1(\tau_1)<\gamma_1(\tau_2)<\gamma_1(\tau_3)$, contradicting the property $\gamma(0)=\gamma(\tau_1)\in(\gamma(\tau_2),\gamma(\tau_3))$ by assumption~\eqref{xitau14bis}. The proof of Lemma~\ref{oneleftbis} is thereby complete.\hfill$\Box$\break

Lemma~\ref{logbis} below is the analogue of Lemma~\ref{log} above, but it is concerned with the logarithmic growth of the argument $\varphi$ of $\nabla u$ along the trajectories of the flow, that is, along the streamlines. Together with Lemmas~\ref{foliation} and~\ref{log}, it will easily lead to the conclusion of Proposition~\ref{growth} on the logarithmic growth of the oscillations of the arguments $\phi$ and $\varphi$ of $w$ and $\nabla u$ in large balls.\par
To state Lemma~\ref{logbis}, let us first introduce a few auxiliary constants. Denote
\beq\label{defomega}
\omega:=\frac{\eta^2}{2}\in\Big(0,\frac12\Big]
\eeq
and 
\beq\label{defC5}
C_2(\eta):=\max\Big(\sup_{t\in[1,+\infty)}h_2(t),288\pi(2\eta^{-4}+1)\Big),\ \hbox{ where }h_2(t)=384\pi\times\frac{\ln_{1+\omega}(4t\eta^{-2})+1}{\ln(3+t)}.
\eeq
Notice immediately that
\beq\label{C5bis}
C_2(\eta)\ge\lim_{t\to+\infty}h_2(t)=\frac{384\pi}{\ln(1+\omega)}\ge\frac{384\pi}{\ln(3/2)}>192\pi>\frac{\pi}{4}.
\eeq
With this constant $C_2(\eta)>0$, the following estimate holds for any streamline of the flow.

\begin{lem}\label{logbis}
For any $\eta\in(0,1]$, for any $C^2(\R^2)$ solution $w$ of~\eqref{1} satisfying~\eqref{etaw} and for any $x\in\R^2$, there holds
$$\mathop{\osc}_{B(x,1)}\varphi=\mathop{\osc}_{B(x,1)}\phi<\frac{\pi}{4}\ \Longrightarrow\ \Big(\forall\,y\in\Gamma_x,\ \ |\varphi(x)-\varphi(y)|\le C_2(\eta)\,\ln\!\big(3+|x-y|\big)\Big).$$
\end{lem}

The scheme of the proof is similar to that of Lemma~\ref{log}. However, since $u$ is constant along any streamline, one cannot conclude directly as in Lemma~\ref{log} that if a streamline turns many times around some points then the stream function $u$ would become large. Nevertheless, roughly speaking, if on some streamline an arc $\gamma([\alpha,\beta])$ contains either two left or two right or two double arcs and is such that $\gamma(\alpha)$ and $\gamma(\beta)$ are not too far, then $\gamma([\alpha,\beta])$ will almost surround a domain~$\Omega$ containing a long trajectory of the gradient flow. This will lead to a contradiction, since the oscillations of $u$ on the boundary of that domain~$\Omega$ will be small, while the oscillation of $u$ along the long trajectory of the gradient flow will be large.\hfill\break

\noindent{\bf{Proof.}} Assume by contradiction that the conclusion of Lemma~\ref{logbis} does not hold for some $C^2(\R^2)$ solution $w$ of~\eqref{1} satisfying~\eqref{etaw} with $0<\eta\le1$. Then, there exist $x\in\R^2$ and $y\in\Gamma_x$ such that ${\osc}_{B(x,1)}\varphi<\pi/4$ and
\beq\label{absurde}
|\varphi(x)-\varphi(y)|>C_2(\eta)\,\ln\!\big(3+|x-y|\big).
\eeq
Let $\gamma$ be a parametrization of $\Gamma_x$ solving $\dot\gamma(t)=w(\gamma(t))$ for $t\in\R$, and $a,b\in\R$ be such that
$$x=\gamma(a)\ \hbox{ and }\ y=\gamma(b).$$
One has $|\varphi(x)-\varphi(y)|>C_2(\eta)>\pi/4$ (by~\eqref{C5bis}) and the property ${\osc}_{B(x,1)}\varphi<\pi/4$ implies that
\beq\label{gammaab}
|\gamma(a)-\gamma(b)|=|x-y|\ge1.
\eeq
Let us assume here that $a<b$ (the case $a>b$ can be handled similarly). Let $\theta$ be the continuous argument of $\dot\gamma$, defined as in~\eqref{deftheta} with the embedding $\xi:=\gamma$. Since $\dot\gamma(t)=w(\gamma(t))$ in $\R$, it follows from~\eqref{defvarphi}, from the continuity of $\theta$ and $\varphi(\gamma(\cdot))$ in $\R$ and from the definition $w=\nabla^{\perp}u$ that there is an integer $q\in\Z$ such that
$$\theta(t)=\varphi(\gamma(t))+\frac{\pi}{2}+2\pi q\ \hbox{ for all }t\in\R.$$
In particular, $\theta(a)-\theta(b)=\varphi(\gamma(a))-\varphi(\gamma(b))=\varphi(x)-\varphi(y)$, hence
\beq\label{thetagamma}
|\theta(a)-\theta(b)|>C_2(\eta)\,\ln\!\big(3+|x-y|\big)=C_2(\eta)\,\ln\big(3+|\gamma(a)-\gamma(b)|\big)>C_2(\eta)>192\pi
\eeq
by~\eqref{C5bis} and~\eqref{absurde}.\par
With the same notations as in Section~\ref{sec31}, let $N_l$, $N_r$ and $N_d$ be the (finite) numbers of left, right and double arcs contained in the curve $\gamma([a,b])$, relatively to the segment $[\gamma(a),\gamma(b)]$. As for~\eqref{theta16} and~\eqref{theta96}, it follows from the inequalities $|\theta(a)-\theta(b)|>192\pi>96\pi>8\pi$ and from Lemma~\ref{arcs} that
$$\max\big(N_l,N_r,N_d\big)>\frac{|\theta(a)-\theta(b)|}{96\pi}.$$
As in the proof of Lemma~\ref{log}, three cases can then occur.\par
\vskip 0.2cm
{\it Case 1: $N_l>|\theta(a)-\theta(b)|/(96\pi)$.}  Notice immediately that~\eqref{thetagamma} yields $N_l>2$, that is, $N_l\ge3$. By definition of $N_l$ and of a left arc, there are some real numbers
$$a<\rho_1<\rho'_1\le\rho_2<\rho'_2\le\cdots\le\rho_{N_l}<\rho'_{N_l}\le b$$
such that $\gamma([\rho_i,\rho'_i])$ is a left arc for every $1\le i\le N_l$.\par
Define
\beq\label{defmbis}
m=\big[\ln_{1+\omega}(4|\gamma(a)-\gamma(b)|\eta^{-2})\big]+2,
\eeq
where $\omega>0$ is defined in~\eqref{defomega}. Property~\eqref{gammaab} together with $0<\eta\le1$ implies that $m$ is an integer such that $m\ge2$. Divide now the segment $[\gamma(a),\gamma(b)]$ into $m$ segments $[x_p,y_p]$ (for $p=1,\cdots,m$) defined by
$$\left\{\baa{ll}
\displaystyle x_1=\gamma(a)+\frac{\gamma(b)-\gamma(a)}{1+\omega},\ y_1=\gamma(b), & \ \ \displaystyle x_2=\gamma(a)+\frac{\gamma(b)-\gamma(a)}{(1+\omega)^2},\ y_2=x_1,\vspace{3pt}\\
\cdots & \\
\displaystyle x_{m-1}=\gamma(a)+\frac{\gamma(b)-\gamma(a)}{(1+\omega)^{m-1}},\ y_{m-1}=x_{m-2}, & \ \ \displaystyle x_m=\gamma(a)\hbox{ and }y_m=x_{m-1}.\eaa\right.$$
It is immediate to see that the sets $(x_p,y_p]$ (with $p\in\{1,\ldots,m\}$) form a partition of $(\gamma(a),\gamma(b)]$, that
\beq\label{xp-yp}
|x_p-y_p|=\omega(1+\omega)^{-p}|\gamma(a)-\gamma(b)|\ \hbox{ for each }1\le p\le m-1
\eeq
and that
$$|\gamma(a)-y_m|=|x_m-y_m|=(1+\omega)^{-(m-1)}|\gamma(a)-\gamma(b)|<\frac{\eta^2}{4}$$
by~\eqref{defmbis}.\par
For each $1\le i\le N_l$, the arc $\gamma([\rho_i,\rho'_i])$ is a left arc, hence either $\gamma(\rho_i)$ or $\gamma(\rho'_i)$ belongs to $(\gamma(a),\gamma(b)]=\bigcup_{1\le p\le m}(x_p,y_p]$. Therefore, by calling, for $1\le p\le m$,
$$J_p:=\big\{i\in\{1,\ldots,N_l\};\ \gamma(\rho_i)\in(x_p,y_p]\hbox{ or }\gamma(\rho'_i)\in(x_p,y_p]\big\},$$
it follows that
\beq\label{Jpbis}
\bigcup_{1\le p\le m}J_p=\{1,\ldots,N_l\}.
\eeq
We now claim that
\beq\label{leftarcbis}
\#J_m\le 1,
\eeq
that is, there is at most one left arc intersecting $(x_m,y_m]=(\gamma(a),y_m]$. Assume by contradiction that there are two integers $i$ and $j$ with $1\le i<j\le N_l$ and such that either $\gamma(\rho_i)$ or $\gamma(\rho'_i)$ belongs to $(x_m,y_m]=(\gamma(a),y_m]$ and either $\gamma(\rho_j)$ or $\gamma(\rho'_j)$ belongs to $(\gamma(a),y_m]$. Since $|\gamma(a)-y_m|<\eta^2/4$, one infers in particular that
$$\min\big(|\gamma(\rho_j)-\gamma(a)|,|\gamma(\rho'_j)-\gamma(a)|\big)<\frac{\eta^2}{4}.$$
Furthermore, $\gamma(a)\in(\gamma(\rho_i),\gamma(\rho'_i))$ since $\gamma([\rho_i,\rho'_i])$ is a left arc. Since $\gamma(a)=x$ and ${\osc}_{B(x,1)}\varphi<\pi/4$, Lemma~\ref{oneleftbis} applied with
$$(\tau_1,\tau_2,\tau_3)=(a,\rho_i,\rho'_i)$$
and $\tau_4=\rho_j$ if $|\gamma(\rho_j)-\gamma(a)|<\eta^2/4$ (resp. $\tau_4=\rho'_j$ if $|\gamma(\rho'_j)-\gamma(a)|<\eta^2/4$) leads to a contradiction. Therefore, the claim~\eqref{leftarcbis} follows.\par
Putting together~\eqref{Jpbis} and~\eqref{leftarcbis}, it follows that there is $p\in\{1,\ldots,m-1\}$ such that
\beq\label{defnbis}
\#J_p\ge\frac{N_l-1}{m-1}
\eeq
(we also remind that $N_l\ge 3$ and $m\ge 2$, whence $\#J_p$ is a positive integer). We claim that $\#J_p\ge 3$. Indeed, notice first that
$$N_l-1>\frac{|\theta(a)-\theta(b)|}{96\pi}-1>\frac{|\theta(a)-\theta(b)|}{192\pi}>\frac{C_2(\eta)\,\ln(3+|\gamma(a)-\gamma(b)|)}{192\pi}$$
by~\eqref{thetagamma}. From~\eqref{defnbis} and the definition of $m$ in~\eqref{defmbis}, one infers that
$$\#J_p\ge\frac{N_l-1}{m-1}>\frac{C_2(\eta)}{192\pi}\times\frac{\ln(3+|\gamma(a)-\gamma(b)|)}{\ln_{1+\omega}(4|\gamma(a)-\gamma(b)|\eta^{-2})+1}.$$
Therefore, $\#J_p>2C_2(\eta)/h_2(|\gamma(a)-\gamma(b)|)\ge2$ by~\eqref{defC5} and~\eqref{gammaab}. Thus, $\#J_p>2$, that is, $\#J_p\ge3$.\par
As a consequence, there are three integers
$$i<j<k\ \in J_p,$$
with $a<\rho_i<\rho'_i\le\rho_j<\rho'_j\le\rho_k<\rho'_k\le b$. In particular, since $j\in J_p$, one has either $\gamma(\rho_j)\in(x_,y_p]$ or $\gamma(\rho'_j)\in(x_p,y_p]$. These two cases will be considered in succession.\par
{\it Subcase 1.1:} consider first the case where
$$\gamma(\rho_j)\in(x_p,y_p].$$
Since $\gamma([\rho_j,\rho'_j])$ is a left arc and $\gamma(\rho_j)\in(x_p,y_p]\subset(\gamma(a),\gamma(b)]$, it follows that
$$\gamma(t)\not\in L_{\gamma(a),\gamma(b)}\hbox{ for all }t\in(\rho_j,\rho'_j),\ \hbox{ and }\ \gamma(\rho'_j)\hbox{ is on the left of }\gamma(a)\hbox{ with respect to }\gamma(b)$$
(hence, $\gamma(\rho'_j)\not\in[x_p,y_p]$). On the other hand, either $\gamma(\rho_k)$ or $\gamma(\rho'_k)$ belongs to the segment $[x_p,y_p]$, while $\rho'_j\le\rho_k<\rho'_k$. Consequently, there is
$$\tau\in(\rho'_j,\rho'_k]$$
such that
$$\gamma(\tau)\in[x_p,y_p]\ \hbox{ and }\ \gamma(t)\not\in[x_p,y_p]\hbox{ for all }t\in(\rho_j,\tau),$$
see Figure~7. Thus, $\gamma([\rho_j,\tau])\cap[x_p,y_p]=\{\gamma(\rho_j),\gamma(\tau)\}$. Let $\Omega$ be the non-empty domain surrounded by the closed simple curve $\gamma([\rho_j,\tau])\,\cup\,(\gamma(\rho_j),\gamma(\tau))$ ($\Omega$ is the hatched region in Figure~7) and let $\sigma=\sigma_{\gamma(\rho'_j)}$ be the solution of
\begin{figure}\centering
\includegraphics[scale=0.8]{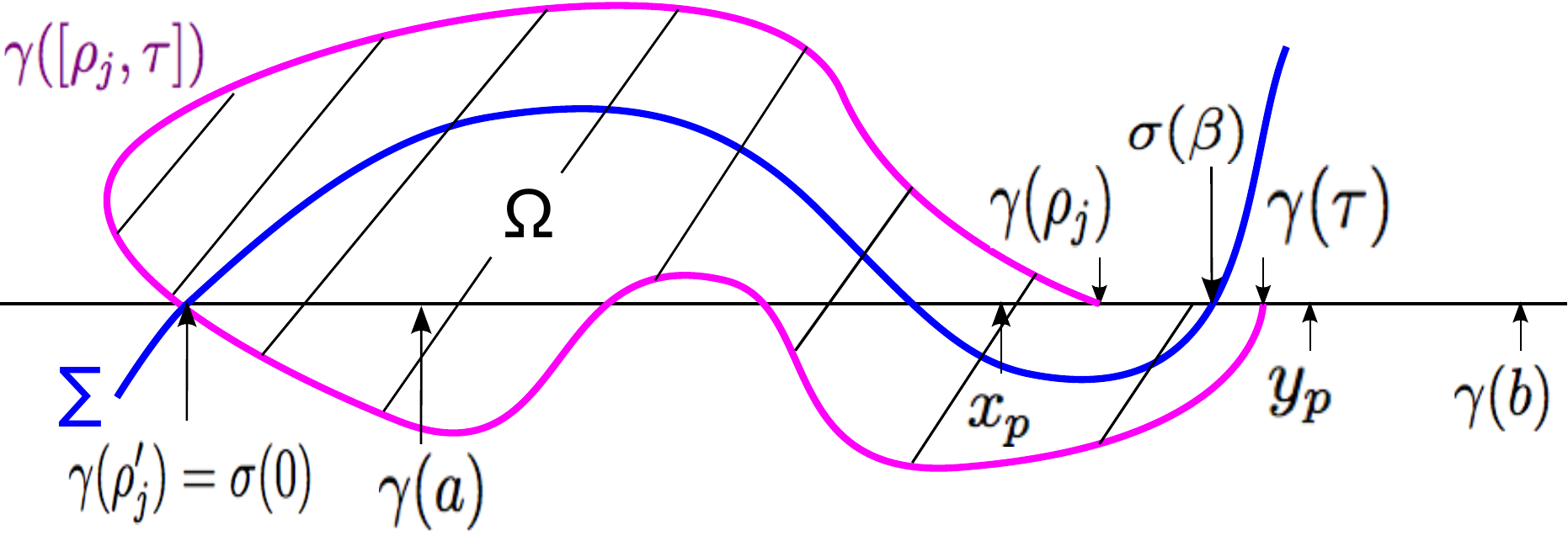}
\caption{The domain $\Omega$ and the curves $\gamma([\rho_j,\tau])$ and $\Sigma$}
\end{figure}
$$\left\{\baa{rcl}
\dot\sigma(t) & = & \nabla u(\sigma(t)),\ \ t\in\R,\vspace{3pt}\\
\sigma(0) & = & \gamma(\rho'_j).\eaa\right.$$
Let $\Sigma=\Sigma_{\gamma(\rho'_j)}=\big\{\sigma(t);\ t\in\R\big\}$ be the trajectory of the gradient flow containing the point $\gamma(\rho'_j)$. Since $\rho_j<\rho'_j<\tau$ and $\Sigma$ is orthogonal to $\Gamma$ at $\gamma(\rho'_j)$ with $\rho'_j\in(\rho_j,\tau)$ and since $|\sigma(t)|\to+\infty$ as $|t|\to+\infty$ by Lemma~\ref{unboundedbis}, it follows that there is a real number $\beta\neq0$ such that
$$\sigma(\beta)\in\partial\Omega\ \hbox{ and }\ \sigma(t)\in\Omega\hbox{ for all }t\in I,$$
where $I=(0,\beta)$ if $\beta>0$ (resp. $I=(\beta,0)$ if $\beta<0$). Since the function $u(\sigma(\cdot))$ is increasing in~$\R$ and $u$ is constant along the curve $\gamma([\rho_j,\tau])\,(\ni\gamma(\rho'_j)=\sigma(0))$, the first exit point $\sigma(\beta)$ must lie on the segment $[\gamma(\rho_j),\gamma(\tau)]$, owing to the definition of $\Omega$. In particular, $\sigma(\beta)\in[x_p,y_p]$. But the points $\sigma(0)=\gamma(\rho'_j)$, $\gamma(a)$, $x_p$, $\sigma(\beta)$ and $y_p$ lie with this order on the line $L_{\gamma(a),\gamma(b)}$. Therefore,
$${\length}(\sigma(I))\ge|\sigma(0)-\sigma(\beta)|\ge|\gamma(a)-x_p|=(1+\omega)^{-p}|\gamma(a)-\gamma(b)|.$$
Thus, Lemma~\ref{lemsigma} yields
\beq\label{inequ1}
|u(\sigma(\beta))-u(\sigma(0))|\ge\eta\times{\length}(\sigma(I))\ge\eta\,(1+\omega)^{-p}|\gamma(a)-\gamma(b)|.
\eeq
But $\sigma(\beta)\in[\gamma(\rho_j),\gamma(\tau)]\subset[x_p,y_p]$, hence~\eqref{etaw} and~\eqref{xp-yp} imply that
\beq\label{inequ2}
|u(\sigma(\beta))-u(\gamma(\rho_j))|\le\eta^{-1}|\sigma(\beta)-\gamma(\rho_j)|\le\eta^{-1}|x_p-y_p|=\eta^{-1}\omega\,(1+\omega)^{-p}|\gamma(a)-\gamma(b)|.
\eeq
Since $\sigma(0)=\gamma(\rho'_j)$ and $u(\sigma(0))=u(\gamma(\rho'_j))=u(\gamma(\rho_j))$, we finally get from~\eqref{inequ1} and~\eqref{inequ2} that
$$\eta\,(1+\omega)^{-p}|\gamma(a)-\gamma(b)|\le\eta^{-1}\omega\,(1+\omega)^{-p}|\gamma(a)-\gamma(b)|,$$
hence $(0<)\,\eta\le\eta^{-1}\omega$, contradicting the definition~\eqref{defomega} of $\omega$.\par
{\it Subcase 1.2:} consider now the case where
$$\gamma(\rho'_j)\in(x_p,y_p].$$
In this case, $\gamma(\rho_j)$ is on the left of $\gamma(a)$ with respect to $\gamma(b)$, while either $\gamma(\rho_i)$ or $\gamma(\rho'_i)$ belongs to $(x_p,y_p]\,(\subset(\gamma(a),\gamma(b)])$. Similarly as in Subcase~1.1, there is $\nu\in[\rho_i,\rho_j)$ such that
$$\gamma(\nu)\in[x_p,y_p]\ \hbox{ and }\ \gamma(t)\not\in[x_p,y_p]\hbox{ for all }t\in(\nu,\rho'_j).$$
We then consider the non-empty domain $\tilde{\Omega}$ surrounded by the closed simple curve $\gamma([\nu,\rho'_j])\cup(\gamma(\nu),\gamma(\rho'_j))$, and the trajectory $\tilde\Sigma$ of the gradient flow containing $\gamma(\rho_j)$. Then $\tilde{\Sigma}\cap\tilde{\Omega}$ contains a connected component having as end points the point $\gamma(\rho_j)$ and a point $z$ belonging to the segment $[\gamma(\nu),\gamma(\rho'_j)]\,(\subset[x_p,y_p])$. As in Subcase~1.1, we similarly get a contradiction by estimating $|u(z)-u(\gamma(\rho_j))|$ from below and above.\par
\vskip 0.2cm
{\it Case 2: $N_r>|\theta(a)-\theta(b)|/(96\pi)$.} The study of this case is similar to that of Case~1 (one especially uses Lemma~\ref{oneleftbis} with the second alternative in~\eqref{tau14bis} to show that the number of right arcs meeting a small portion of $[\gamma(a),\gamma(b)]$ around the point $\gamma(b)$ is at most~$1$), and one similarly gets a contradiction with the definition of $\omega$.
\vskip 0.2cm
{\it Case 3: $N_d>|\theta(a)-\theta(b)|/(96\pi)$.} As in Case~1, one then has $N_d\ge3$. By definition of $N_d$ and of a double arc, there are some real numbers
$$a<\mu_1<\mu'_1\le\mu_2<\mu'_2\le\cdots\le\mu_{N_d}<\mu'_{N_d}\le b$$
such that $\gamma([\mu_i,\mu'_i])$ is a left arc for every $1\le i\le N_d$. For any such $i$, denote $\mu^-_i$ the one of the two real numbers $\mu_i$ and $\mu'_i$ such that $\gamma(\mu^-_i)$ lies on the left of $\gamma(a)$ with respect to $\gamma(b)$, and let $\mu^+_i$ be the other one, that is, $\gamma(b)\in(\gamma(a),\gamma(\mu^+_i))$. We will argue a bit as in Case~1 above, by considering two subcases according to the value of the leftmost position among the points $\gamma(\mu^-_i)$ (notice that $(\gamma(\mu^-_i)-\gamma(a))\cdot(\gamma(b)-\gamma(a))<0$ for all $1\le i\le N_d$, by definition of $\mu^-_i$).\par
{\it Subcase 3.1:} consider first the case where
$$\min_{1\le i\le N_d}\,(\gamma(\mu^-_i)-\gamma(a))\cdot(\gamma(b)-\gamma(a))\le-2\,\eta^{-2}|\gamma(a)-\gamma(b)|^2.$$
In other words, there is $i\in\{1,\ldots,N_d\}$ such that
\beq\label{eta-2}
|\gamma(\mu^-_i)-\gamma(a)|\ge2\,\eta^{-2}|\gamma(a)-\gamma(b)|.
\eeq
Since $\gamma([\mu_i,\mu'_i])$ is a double arc, one has $[\gamma(a),\gamma(b)]\subset(\gamma(\mu_i),\gamma(\mu'_i))$ and $\gamma([\mu_i,\mu'_i])\cap[\gamma(a),\gamma(b)]=\emptyset$. Therefore, there are some real numbers $\lambda$ and $\lambda'$ such that
\begin{figure}\centering
\includegraphics[scale=0.8]{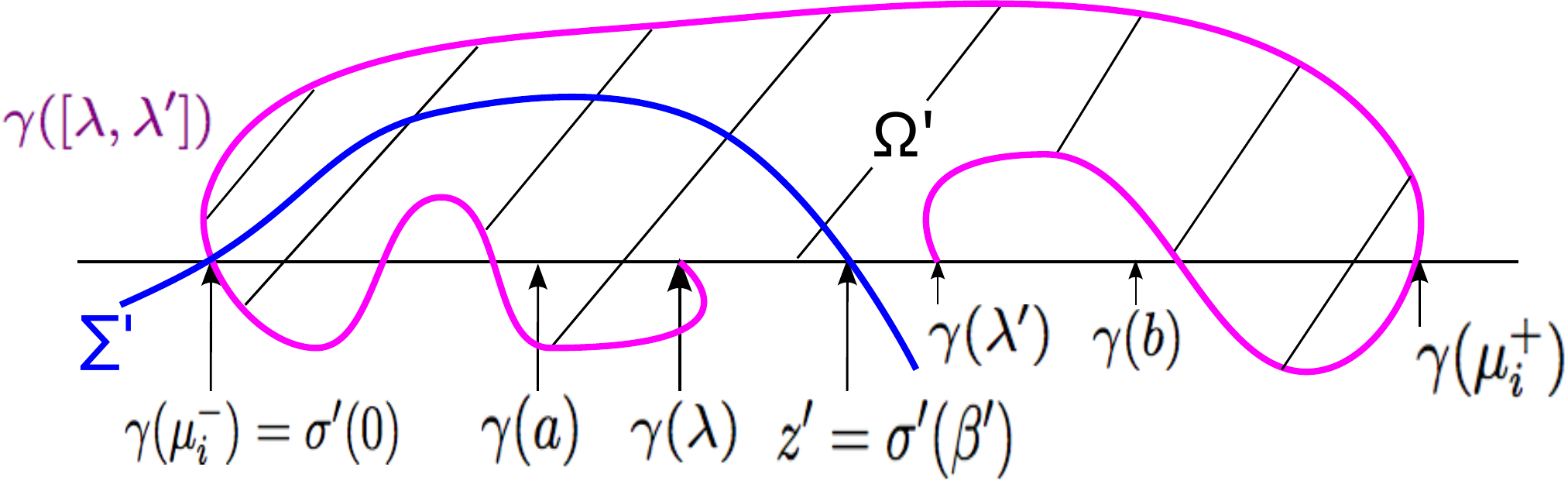}
\caption{The domain $\Omega'$ and the curves $\gamma([\lambda,\lambda'])$ and $\Sigma'$}
\end{figure}
$$a\le\lambda<\mu_i<\mu'_i<\lambda'\le b\ \hbox{ and }\ \gamma([\lambda,\lambda'])\cap[\gamma(a),\gamma(b)]=\{\gamma(\lambda),\gamma(\lambda')\},$$
see Figure~8 (in case $\mu^-_i=\mu_i$ and $\mu^+_i=\mu'_i$). Let $\Omega'$ be the non-empty domain surrounded by the closed simple curve $\gamma([\lambda,\lambda'])\cup(\gamma(\lambda),\gamma(\lambda'))$, let $\Sigma'$ be the trajectory of the gradient flow containing $\gamma(\mu_i^-)$ (which lies on the left of $\gamma(a)$ with respect to $\gamma(b)$), and let $\sigma'=\sigma_{\gamma(\mu^-_i)}$ be the solution of~\eqref{xix} with $\sigma'(0)=\gamma(\mu^-_i)$. Then $\Sigma'\cap\Omega'$ contains a connected component having as end points the point $\gamma(\mu^-_i)$ and a point $z'=\sigma'(\beta')\in[\gamma(\lambda),\gamma(\lambda')]\,(\subset[\gamma(a),\gamma(b)])$, for some $\beta'\neq 0$. As in Case~1 above, using Lemma~\ref{lemsigma} and the fact that $\gamma(\mu^-_i)$ lies on left of $\gamma(a)$ with respect to $\gamma(b)$, one infers from~\eqref{eta-2} that
$$|u(\sigma'(\beta'))-u(\sigma'(0))|\ge\eta\,|\sigma'(\beta')-\sigma'(0)|=\eta\,|\sigma'(\beta')-\gamma(\mu^-_i)|\ge\eta\,|\gamma(a)-\gamma(\mu^-_i)|\ge2\,\eta^{-1}|\gamma(a)-\gamma(b)|.$$
But $\sigma'(\beta')\in[\gamma(a),\gamma(b)]$, and thus
$$|u(\sigma'(\beta'))-u(\gamma(a))|\le\eta^{-1}|\gamma(a)-\gamma(b)|$$
by~\eqref{etaw}. Since $u(\sigma'(0))=u(\gamma(\mu^-_i))=u(\gamma(a))$, one concludes from the previous two displayed inequalities that $2\,\eta^{-1}|\gamma(a)-\gamma(b)|\le\eta^{-1}|\gamma(a)-\gamma(b)|$, which is impossible since both $\eta$ and $|\gamma(a)-\gamma(b)|$ are positive.\par
{\it Subcase 3.2:} consider now the case where
$$\min_{1\le i\le N_d}\,(\gamma(\mu^-_i)-\gamma(a))\cdot(\gamma(b)-\gamma(a))>-2\,\eta^{-2}|\gamma(a)-\gamma(b)|^2.$$
Denote
$$\zeta=\gamma(a)-2\,\eta^{-2}(\gamma(b)-\gamma(a)).$$
It follows that $\gamma(\mu^-_i)\in(\zeta,\gamma(a)]$ for all $1\le i\le N_d$. Denote now
\beq\label{defN}
N=\big[2\,\eta^{-4}+1\big]\,(\ge 3)
\eeq
and divide the segment $[\zeta,\gamma(a)]$ into $N$ subsegments
$$[x_p,y_p]=\Big[\zeta+\frac{p-1}{N}\,(\gamma(a)-\zeta),\zeta+\frac{p}{N}\,(\gamma(a)-\zeta)\Big]$$
(for $1\le p\le N$) of the same length
\beq\label{eta-2bis}
|x_p-y_p|=\frac{|\gamma(a)-\zeta|}{N}=\frac{2\,\eta^{-2}|\gamma(b)-\gamma(a)|}{N}.
\eeq
Remember that $N_d>|\theta(a)-\theta(b)|/(96\pi)$ in our studied Case~3, and that
$$|\theta(a)-\theta(b)|>C_2(\eta)\ge 288\pi(2\eta^{-4}+1)\ge 288\,\pi\,N$$
by~\eqref{defC5},~\eqref{thetagamma} and~\eqref{defN}. Therefore, $N_d>3N$. In particular, there are $p\in\{1,\ldots,N\}$ and three integers $i<j<k$ in $\{1,\ldots,N_d\}$ such that
$$\gamma(\mu^-_i),\,\gamma(\mu^-_j),\, \gamma(\mu^-_k)\,\in\,[x_p,y_p].$$\par
If $\mu^-_j=\mu_j$, then $\mu^-_j=\mu_j<\mu'_j\le\mu_k<\mu'_k$, while either $\gamma(\mu_k)$ or $\gamma(\mu'_k)$ belongs to $[x_p,y_p]$. Since $\gamma(\mu_j)=\gamma(\mu^-_j)\in[x_p,y_p]$ (lying on the left of $\gamma(a)$ with respect to $\gamma(b)$) and since $\gamma([\mu_j,\mu'_j])$ is a double arc, there is then
$$\varpi\in(\mu'_j,\mu'_k]$$
such that
\begin{figure}\centering
\includegraphics[scale=0.8]{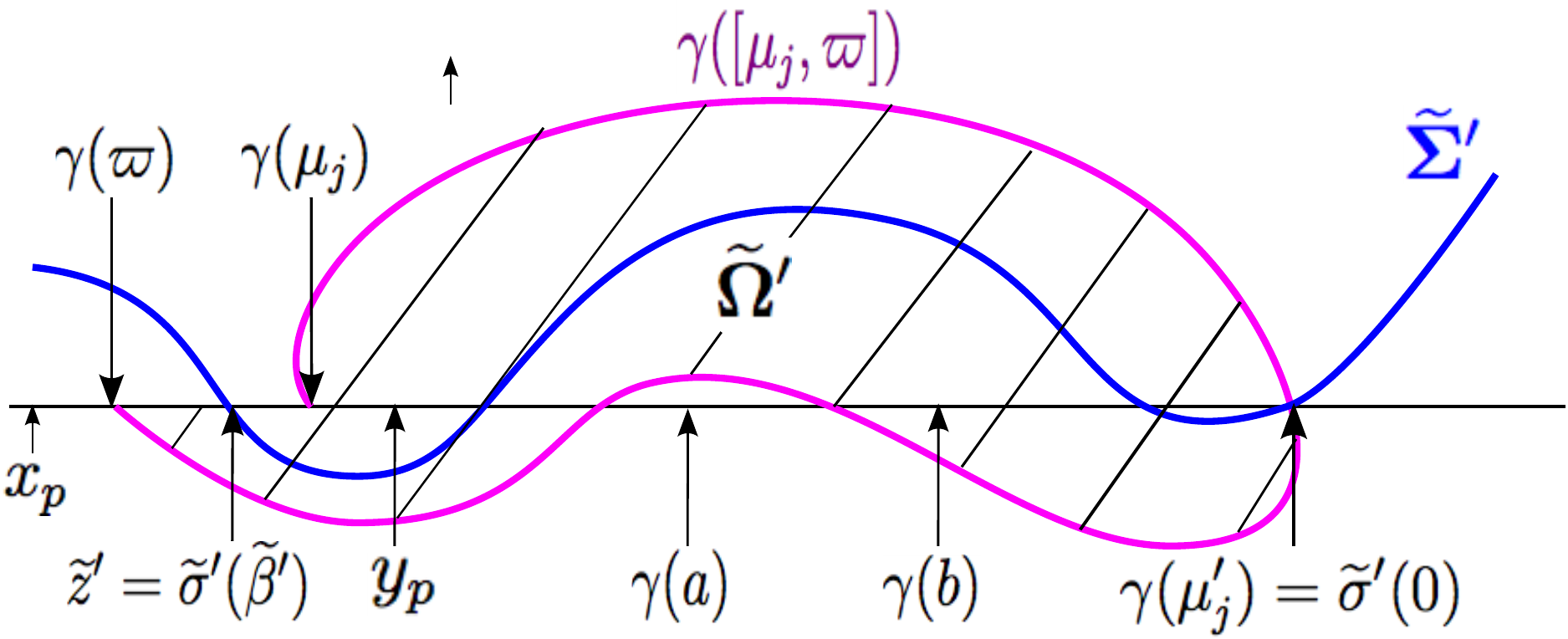}
\caption{The domain $\tilde{\Omega}'$ and the curves $\gamma([\mu_j,\varpi])$ and $\tilde{\Sigma}'$}
\end{figure}
$$\gamma(\varpi)\in[x_p,y_p]\ \hbox{ and }\ \gamma(t)\not\in[x_p,y_p]\hbox{ for all }t\in(\mu_j,\varpi),$$
see Figure~9. Thus, $\gamma([\mu_j,\varpi])\cap[x_p,y_p]=\{\gamma(\mu_j),\gamma(\varpi)\}$. Let $\tilde{\Omega}'$ be the non-empty domain surrounded by the closed simple curve $\gamma([\mu_j,\varpi])\,\cup\,(\gamma(\mu_j),\gamma(\varpi))$ and let $\tilde{\sigma}'=\sigma_{\gamma(\mu'_j)}$ be the solution of~\eqref{xix} with initial condition $\tilde{\sigma}'(0)=\gamma(\mu'_j)$. Let $\tilde{\Sigma}'=\Sigma_{\gamma(\mu'_j)}$ be the trajectory of the gradient flow containing the point $\gamma(\mu'_j)$. There is then a real number $\tilde{\beta}'\neq0$ such that $\tilde{\Sigma}'\cap\tilde{\Omega}'$ contains a connected component having as end points the point $\tilde{\sigma}'(0)=\gamma(\mu'_j)$ and the point $\tilde{z}'=\tilde{\sigma}'(\tilde{\beta}')\in[\gamma(\mu_j),\gamma(\varpi)|\subset[x_p,y_p]$. But the points $x_p$, $\tilde{z}'=\tilde{\sigma}'(\tilde{\beta}')$, $y_p$, $\gamma(a)$, $\gamma(b)$ and $\gamma(\mu'_j)=\tilde{\sigma}'(0)$ lie in this order on the line $L_{\gamma(a),\gamma(b)}$. Therefore, Lemma~\ref{lemsigma} yields
$$|u(\tilde{\sigma}'(\tilde{\beta}'))-u(\tilde{\sigma}'(0))|\ge\eta\,|\tilde{\sigma}'(\tilde{\beta}')-\tilde{\sigma}'(0)|\ge\eta\,|\gamma(a)-\gamma(b)|.$$
But $\tilde{\sigma}'(\tilde{\beta}')$ and $\gamma(\mu_j)$ belong to the segment $[x_p,y_p]$. Hence,
$$|u(\tilde{\sigma}'(\tilde{\beta}'))-u(\gamma(\mu_j))|\le\eta^{-1}|\tilde{\sigma}'(\tilde{\beta}')-\gamma(\mu_j)|\le\eta^{-1}|x_p-y_p|=\frac{2\,\eta^{-3}|\gamma(a)-\gamma(b)|}{N}$$
by~\eqref{etaw} and~\eqref{eta-2bis}. Since $u(\tilde{\sigma}'(0))=u(\gamma(\mu'_j))=u(\gamma(\mu_j))$ and $|\gamma(a)-\gamma(b)|>0$, the previous two inequalities imply that $\eta\le2\eta^{-3}/N$, contradicting~\eqref{defN}.\par
Finally, if $\mu^-_j=\mu'_j$, then $\mu_i<\mu'_i\le\mu_j<\mu'_j=\mu^-_j$, while either $\gamma(\mu_i)$ or $\gamma(\mu'_i)$ belongs to $[x_p,y_p]$. Since $\gamma(\mu'_j)=\gamma(\mu^-_j)\in[x_p,y_p]$ (lying on the left of $\gamma(a)$ with respect to $\gamma(b)$) and since $\gamma([\mu_j,\mu'_j])$ is a double arc, there is then $\varpi'\in[\mu_i,\mu_j)$ such that $\gamma(\varpi')\in[x_p,y_p]$ and $\gamma(t)\not\in[x_p,y_p]$ for all $t\in(\varpi',\mu'_j)$. Thus, $\gamma([\varpi',\mu'_j])\cap[x_p,y_p]=\{\gamma(\varpi'),\gamma(\mu'_j)\}$ and one gets a contradiction with similar arguments as in the previous paragraph.\par
To sum up, all possible cases have been considered. They all lead to a contradiction. As a conclusion,~\eqref{absurde} cannot hold and the proof of Lemma~\ref{logbis} is thereby complete.\hfill$\Box$


\subsection{End of the proof of Proposition~\ref{growth}}\label{sec34}

With Lemmas~\ref{log} and~\ref{logbis} in hand, the proof of Proposition~\ref{growth} follows easily. To do so, consider any $\eta\in(0,1]$, any $C^2(\R^2)$ solution $w$ of~\eqref{1} satisfying~\eqref{etaw}, any real number $R\ge2$, and assume that ${\osc}_{B(x,1)}\varphi_w={\osc}_{B(x,1)}\phi_w<\pi/4$ for all $x\in B(0,R)$, with the general notations of Section~\ref{sec22}. Let $\Sigma=\Sigma_0$ be the trajectory of the gradient flow containing the origin~$0=(0,0)$.\par
Let now $x$ be any point in $B(0,R)$. Lemma~\ref{foliation} yields the existence of a point $y_x\in\Sigma$ such that $x\in\Gamma_{y_x}$ (that is, $y_x\in\Gamma_x$). Since ${\osc}_{B(x,1)}\varphi_w<\pi/4$ and ${\osc}_{B(0,1)}\varphi_w<\pi/4<\pi/2$, Lemmas~\ref{log} and~\ref{logbis} imply that
\beq\label{phix}\baa{rcl}
|\varphi_w(x)-\varphi_w(0)| & \le & |\varphi_w(x)-\varphi_w(y_x)|+|\varphi_w(0)-\varphi_w(y_x)|\vspace{3pt}\\
& \le & C_2(\eta)\,\ln(3+|x-y_x|)+C_1(\eta)\,\ln(3+|y_x|).\eaa
\eeq
Furthermore, on the one hand, Lemma~\ref{lemsigma} and the normalization $u(0)=0$ yield
$$|u(y_x)|=|u(y_x)-u(0)|\ge\eta\,|y_x|.$$
On the other hand, $u(y_x)=u(x)$ (since $x\in\Gamma_{y_x}$), hence
$$|u(y_x)|=|u(x)-u(0)|\le\eta^{-1}|x|$$
by~\eqref{etaw}. Therefore, $|y_x|\le\eta^{-2}|x|$. Together with~\eqref{phix} and $|x|<R$, one finally concludes that
\beq\label{oscx0}
|\varphi_w(x)-\varphi_w(0)|\le C_2(\eta)\,\ln(3+|x|+\eta^{-2}|x|)+C_1(\eta)\,\ln(3+\eta^{-2}|x|)\le\frac{C_\eta}{2}\,\ln R,
\eeq
where
$$C_\eta=2\,\sup_{\rho\in[2,+\infty)}\frac{C_2(\eta)\,\ln(3+\rho+\eta^{-2}\rho)+C_1(\eta)\,\ln(3+\eta^{-2}\rho)}{\ln\rho}.$$
Notice that the constant $C_\eta$ is a positive real number depending on $\eta$ only. Since~\eqref{oscx0} holds for every $x\in B(0,R)$, one concludes that
$$\mathop{\osc}_{B(0,R)}\phi_w=\mathop{\osc}_{B(0,R)}\varphi_w\le C_\eta\ln R.$$
The proof of Proposition~\ref{growth} is thereby complete.\hfill$\Box$


\end{document}